\documentclass{article}
\usepackage{amsmath,amsthm,amsfonts,amssymb}
\usepackage[all]{xypic}
\usepackage[final]{graphicx}

\usepackage{stmaryrd} 

\usepackage{hyperref} 

\usepackage{mathrsfs} 

\newcommand{\normal}{{\trianglelefteqslant}}

\newcommand{\setR}{\mathbb R}
\newcommand{\setC}{\mathbb C}
\newcommand{\setH}{{\mathbb H}}

\newcommand{\g}{\mathfrak g}

\newcommand\p{{\partial}}
\newcommand\G{{\Gamma}}
\newcommand\OC{{\overline C}}
\newcommand\U{{\mathcal U}}

\newcommand{\C}[1]{{\cal#1}} 

\theoremstyle{plain}
\newtheorem{prop}{Proposition}[section]
\newtheorem{thm}{Theorem}[section]
\newtheorem{lem}{Lemma}[section]

\theoremstyle{definition}
\newtheorem{rem}{Remark}[section]
\newtheorem{example}{Example}[section]
\newtheorem{definition}{Definition}[section]

\numberwithin{equation}{section}

\begin{document}

\title{Graph complexes in deformation quantization}

\author{Domenico Fiorenza and Lucian M. Ionescu}

\maketitle

\begin{abstract}
Kontsevich's formality theorem and the consequent star-product formula 
rely on the construction of an $L_\infty$-morphism between 
the DGLA of polyvector fields and the DGLA of polydifferential 
operators. 
This construction uses a version of graphical calculus. 
In this article we present the details of this graphical calculus 
with emphasis on its algebraic features.
It is a morphism of differential graded Lie algebras between 
the Kontsevich DGLA of admissible graphs and the 
Chevalley-Eilenberg DGLA of linear homomorphisms
between polyvector fields and polydifferential operators. 
\par
Kontsevich's proof of the formality morphism is reexamined in this 
light and an algebraic framework for discussing the tree-level 
reduction of Kontsevich's star-product is described.
\end{abstract}

\tableofcontents

\section{Introduction}

In the breakthrough paper \cite{K1} solving the long-standing
problem of  deformation quantization, 
Maxim Kontsevich
makes use of a perturbative approach in writing his
solution: the formality
$L_\infty$-quasi-isomorphism. The ideological background is
that of a string theory and, as it was explained in
\cite{CF},   the formality morphism is really the
perturbative expansion of the partition function of a sigma
model on the unit disk with the 2-point function yielding
the star-product. From a mathematical point of view,
Kontsevich's construction can be seen as a graphical
calculus for derivations, providing another example of the
ubiquity of graphical calculi in contemporary mathematics 
(e.g. Feynman diagrams, Turaev's tangles, knot  and
3-manifold invariants, TQFTs, operads and PROPs etc.). The
analogy between Kontsevich's and Feynman's graphical
calculus was further pointed out in
\cite{pqft}.
\par
In this article we focus on the graphical calculus for 
derivations
from  the mathematical point of view, insisting on the fact 
that the Kontsevich's rule implements graphically an action,
which extends the basic rule associating to a colored arrow 
the corresponding evaluation 
$$
\begin{xy}
,(-3,3);(3,-3)*{\bullet}**\crv{(-2,-3)}?>*\dir{>}
,(-3.1,-1)*{\scriptstyle{\varphi}},(3,-4.5)*{\scriptstyle{x}}
\end{xy}
\,\,\mapsto \varphi(x)
$$
by incorporating the Leibniz rule relative to the natural tensor
product given by disjoint union:
$$
\begin{xy}
,(-3,3);(3,-3)**\crv{(-2,-3)}?>*\dir{>}
,(6,-3)*\cir(3,3){},
,(-3.1,-1)*{\scriptstyle{\varphi}},(4.5,-4.3)*{\scriptscriptstyle{x}}
,(7.5,-4.3)*{\scriptscriptstyle{y}}
,(4.5,-3)*{\bullet},(7.5,-3)*{\bullet}
\end{xy}
=
\begin{xy}
,(-3,3);(3,-3)*{\bullet}**\crv{(-2,-3)}?>*\dir{>}
,(-3.1,-1)*{\scriptstyle{\varphi}},(3,-4.5)*{\scriptstyle{x}}
,(7,-3)*{\bullet},(7,-4.5)*{\scriptstyle{y}}
\end{xy}
+
\begin{xy}
,(-3,3);(7,-3)*{\bullet}**\crv{(4,3)}?>*\dir{>}
,(5,1.5)*{\scriptstyle{\varphi}},(7,-4.5)*{\scriptstyle{y}}
,(3,-3)*{\bullet},(3,-4.5)*{\scriptstyle{x}}
\end{xy}
\,\,\mapsto \varphi(x)\,y+x\,\varphi(y).
$$

The main result is Theorem \ref{T:DGLA} stating that Kontsevich's
 graphical calculus
is a DGLA morphism between the differential graded Lie algebra
${\mathscr G}^{\, \bullet}$ of Kontsevich's admissible graphs and the
Chevalley-Eilenberg differential graded Lie algebra
$CE^\bullet(T^\bullet_{\rm poly}, D^\bullet_{\rm poly})$ of linear
homomorphisms between polyvector fields and polydifferential
operators.  From this it follows that a solution of the  Maurer-Cartan
equation in ${\mathscr G}^{\, \bullet}$ induces by push-forward a
solution of the  Maurer-Cartan equation in $CE^\bullet(T^\bullet_{\rm
poly}, D^\bullet_{\rm poly})$; therefore it provides an
$L_\infty$-morphism between $T^\bullet_{\rm poly}$ and
$D^\bullet_{\rm poly}$, which in turn can be used to define a 
star-product quantization of Poisson structures.
The interplay between the DGLA structure on admissible graphs and
its dual dg-Lie-coalgebra structure is mentioned,
pointing out its r\^ole in reducing the problem of solving the
Maurer-Cartan equation on graphs to the problem of finding a system of
``weights''
$W_\Gamma$ satisfying the cocycle equation which is dual to the
Maurer-Cartan equation. 
\par
In the final part of the paper, we show how the DGLA
structure on ${\mathscr G}^{\, \bullet}$ provides a nice 
algebraic framework for studying the tree-level
 approximation of Kontsevich's star product formula. Indeed,
by homotopical transfer of structure, the DGLA structure on
${\mathscr G}^{\, \bullet}$ induces an $L_\infty$-algebra
structure on the subcomplex ${\mathscr F}^{\, \bullet}$ of
forest graphs, and an $L_\infty$-morphism $\U_\infty\colon
{\mathscr F}^{\, \bullet}\to CE^\bullet(T^\bullet_{\rm
poly}, D^\bullet_{\rm poly})$, whose principal part is the
Kontsevich's graphical calculus restricted to forest graphs. 
As a consequence one obtains that it is possible to write 
a star-product formula as a sum over forest graphs; 
this formula is given 
by a Kontsevich-type formula involving only forest graphs 
(the so-called tree-level or semi-classical approximation) 
plus non-linear corrections 
corresponding to non simply-connected graphs in the original
Kontsevich's formula. 
\par
There is a considerable overlap between the first few
sections of this paper and 
other expositions of Kontsevich's paper \cite{CI,Kel,DS,AMM}
etc.,
but we tried to reorganize the presentation to clarify 
the DGLA structure on graphs, 
together with the key feature of Kontsevich's approach,
namely that the graphical calculus $\U$ is a DGLA morphism.

\medskip

\noindent{\it Acknowledgments.} The authors thank Benoit
Dherin, Fabio Gavarini and the Referees from Letters in
Mathematical Physics for very useful comments and suggestions on
the first version of this paper.

\section{The graphical calculus for derivations}\label{S:gc}

An oriented graph can be used to build a
differential operator out of a finite set of polyvector
fields. This construction is quite standard (see, for
instance \cite{K1})
and we will recall only its
computational rules here. To begin with, we need to specify
which graphs we do have in mind.

\subsection{Kontsevich's admissible graphs}\label{SS:kag}

By definition, an \emph{admissible graph} 
with $n$ internal vertices and $m$ boundary vertices is a finite oriented
graph $\Gamma$ endowed with the following additional data/
con\-di\-tions:
\begin{enumerate}
\item the vertices of $\Gamma$ are split into the two totally ordered disjoint
subsets of
\emph{internal} and \emph{boundary} vertices;
\item the internal vertices are numbered from $1$ to $n$;
\item the boundary vertices are numbered from $1$ to $m$;
\item no edge stems from a boundary vertex, i.e., boundary vertices are
sinks;
\item an ordering on the edges stemming from each internal
vertex is given; due to the total order on the internal vertices,
this is equivalent to a total order on the set of all 
edges of $\Gamma$;
\item there are no multiple edges, and no loops, i.e. edges starting and
ending on the same vertex.
\end{enumerate}

A good way of visualizing admissible graphs is as graphs
drawn into a disk, with internal vertices in the interior of
the disk and boundary vertices on the boundary, in such a way
that all their edges are
realized as geodesics for the standard hyperbolic metric of the disk. An example
is given in the following figure, depicting an admissible graph
with
$2$ internal  vertices and $3$
boundary vertices:
 
\vbox{
\[
\begin{xy}
,(0,0)*\cir(10,10){},
,(-7,-7)*{\bullet},(7,-7)*{\bullet}
,(1,-3);(-7,-7)**\crv{(0,-4)} ?>*\dir{>}
,(1,-3);(0,-10)*{\bullet}**\crv{(0,-4)} ?>*\dir{>}
,(3,3);(1,-3)**\crv{(3,0)} ?>*\dir{>}
,(3,3);(7,-7)**\crv{(5,-5)} ?>*\dir{>}
,(3,4.4)*{\scriptstyle{1}}
,(.3,-1.9)*{\scriptstyle{2}}
,(1.8,1.0)*{\scriptscriptstyle{1}}
,(-3,-4)*{\scriptscriptstyle{3}}
,(5.9,-2.5)*{\scriptscriptstyle{2}}
,(-1,-6.5)*{\scriptscriptstyle{4}}
,(-8,-9)*{\scriptstyle{1}}
,(0,-12)*{\scriptstyle{2}}
,(8,-9)*{\scriptstyle{3}}
\end{xy}
\]
\vskip -12 pt
\[
\text{\footnotesize Figure 1: an admissible graph}
\]}

\subsection{Graphical calculus for derivations}\label{S:gcfd}
Let now $x^i$ be coordinates on ${\setR}^d$, and let
$\partial_i$ be the corresponding vector fields.  The space
of polyvector fields on ${\setR}^d$ is 
\[
T^\bullet_{\rm poly}=\bigoplus_{k=-1}^\infty T^k_{\rm
poly}=\bigoplus_{k=-1}^\infty
H^0({\setR}^d;\wedge^{k+1}T{{\setR}^d}).
\] 
The graphical calculus for derivations uses a graph
$\Gamma$ with $n$ internal and $m$ boundary vertices,
together with $n$ polyvector fields,
to build a polydifferential operator acting on $m$ smooth functions defined on
${\setR}^d$. In other words, if we denote by $D_{\rm
poly}^{m-1}$ the space of  polydifferential operators acting
on $m$ smooth functions defined on ${\setR}^d$ 
(\S\ref{SS:Dpoly}; see also \cite{CI}, p.24), 
then the graphical calculus for derivations 
associates to a graph $\Gamma$ with $n$ internal and $m$ boundary vertices
a map
\[
\U_\Gamma\colon  \left(T^\bullet_{\rm
poly}\right)^{\otimes n}\to D_{\rm
poly}^{m-1}.
\]
More precisely, if $\nu_r+1$ is the number of edges stemming
from the $r^{\rm th}$ internal vertex of $\Gamma$, then
$\U_\Gamma$ is a map
\[
\U_\Gamma\colon  \bigotimes_{r=1}^nT^{\nu_r}_{\rm
poly}\to D_{\rm poly}^{m-1}.
\]
which is better expressed in the grading preserving form
\[
\U_\Gamma\colon  \bigotimes_{r=1}^nT^{\nu_r}_{\rm
poly}\to (D_{\rm poly}[v(
\Gamma)-1-e(\Gamma)])^{\sum\nu_r}.
\]
where $e(\Gamma)$ and $v(\Gamma)$ are the number of edges and vertices of
$\Gamma$, respectively, and we have used the identities
\begin{align*}
e(\Gamma)&=n+\sum_{r=1}^n\nu_r;
\\
v(\Gamma)&=m+n.
\end{align*} 
This suggests to introduce a bigrading on admissible graphs as
follows.
\begin{definition}
An {\em admissible graph of type $(p,q)$} is an admissible graph $\Gamma$ such that
\begin{enumerate}
\item $p=\chi(\Gamma)-1=v(\Gamma)-e(\Gamma)-1$;
\item the number of internal vertices of $\Gamma$ is $q$. 
\end{enumerate}
The set of isomorphism classes\footnote{An isomorphism of
admissible graphs is an isomorphism of the underlying oriented graphs
which preserves all the additional structures on them.} of admissible graphs of type
$(p,q)$ will be denoted by the symbol ${\mathcal G}^{p,q}$.
\end{definition}
With these notations, $\U$ can be seen as a map
\[
\U\colon {\mathcal G}^{p,q}\to{\rm Hom}\left((T^\bullet_{\rm
poly})^{\otimes q},D_{\rm
poly}^\bullet[p]\right),
\]
where the homomorphisms on the right hand side are grading preserving. We will
return on this in the following sections. 

Now, let us explicitly describe the rules of graphical calculus. 
Such rules define
\[
\langle\U_\Gamma(\xi_1\otimes\cdots\otimes\xi_n)|
f_0\otimes\cdots\otimes
f_{m-1}\rangle
\]
as a state sum. 
Namely, evaluation of the polydifferential
operator above is computed by
considering the {\em decorated graph} corresponding to our data, 
i.e., the graph
$\Gamma$  with the internal vertices 
labeled by the polyvector fields $\{\xi_r\}_{r=1,\dots,n}$
and  the boundary vertices labeled by the functions
$\{f_j\}_{j=1,\dots,m}$.
Now proceed to a state-sum construction;
the edges stemming from each internal vertex
$v$ of $\Gamma$ are labeled by tensor indices $i_{k}$
with $k$ given by the numbering induced by the total order
on the set of edges. 
We refer to this operation as ``coloring'' the
graph.

The value of such a
{\em colored graph} is defined as the product over the
contributions of the vertices of $\Gamma$, where the value
of a decorated vertex $v$ is
\[
\underbrace{\overbrace{
\begin{xy}
,(0,4.5)*{\phantom{m}}
,(-5,4);(0,0)**\dir{-} ?>*\dir{>}
,(-2,4);(0,0)**\dir{-} ?>*\dir{>}
,(.5,3)*{\scriptstyle{\dots}}
,(5,4);(0,0)**\dir{-} ?>*\dir{>}
,(0,0);(-5,-4)**\dir{-}?>*\dir{>}
,(0,0);(5,-4)**\dir{-}?>*\dir{>}
,(0,0);(-2,-4)**\dir{-}?>*\dir{>}
,(.5,-3)*{\scriptstyle{\dots}}
,(3,0)*{\scriptstyle{\xi_v}}
,(0,-4.5)*{\phantom{m}}
\end{xy}}^{{\rm In}_v}}_{{\rm Out}_v}
\mapsto |{\rm Out}_v|!\,
\partial_{{\,\rm In}_v}
\xi_v^{{\rm Out}_v}\,,
\]
for an internal vertex, and
\[
\overbrace{
\begin{xy}
,(0,4.5)*{\phantom{m}}
,(-5,4);(0,0)*{\bullet}**\dir{-} ?>*\dir{>}
,(-2,4);(0,0)*{\bullet}**\dir{-} ?>*\dir{>}
,(.5,3)*{\scriptstyle{\dots}}
,(5,4);(0,0)*{\bullet}**\dir{-} ?>*\dir{>}
,(2,-1)*{\scriptstyle{f_v^{}}}
,(0,-4.5)*{\phantom{m}}
\end{xy}}^{{\rm In}_v}
\mapsto
\partial_{{\,\rm In}_v}
f_v^{}\,,
\]
for a boundary vertex, where ${\rm In}_v$ and ${\rm Out}_v$
are the multi-indices corresponding to the numbering on the
incoming and outgoing edges at the vertex $v$.
Finally, we sum over all repeated indices according to
Einstein's convention.

\begin{rem}
Polydifferential operators can be built by two basic
operations: multiplication of functions, and $m$-vector fields
acting as polyderivations on $m$ functions. These two
basic operations are represented in the graphical calculus
for derivations $\U$ by the most basic graphs. Namely,
multiplication of $m$ functions is represented by the graph
$b_{0,m}$ with no internal vertices and $m$ boundary
vertices:
\[
\langle\U_{b_{0,m}}|f_0^{}\otimes
f_1^{}\otimes\cdots\otimes f_{m-1}^{}\rangle=
\begin{xy}
,(0,0)*\cir(10,10){},
,(-7,-7)*{\bullet},(7,-7)*{\bullet}
,(-4,-9)*{\bullet}
,(-8,-9)*{\scriptstyle{f_0^{}}}
,(-4,-11)*{\scriptstyle{f^{}_1}}
,(1,-8)*{\cdots}
,(9,-9)*{\scriptstyle{f^{}_{m-1}}}
\end{xy}\,=
f_0^{}\cdot f_1^{}\cdots f_{m-1}^{}
\]
It is an element of ${\mathcal G}^{m-1,0}$. 
The identification of $m$-vector fields with
polyderivations is represented by the admissible graph $b_{1,m}$
with one internal vertex, $m$ boundary vertices and $m$
edges (we will call $b_{1,m}$ an {\em $m$-corolla}):
\begin{align*}
\langle\U_{b_{1,m}}(\xi)|f_0^{}\otimes
f_1^{}\otimes&\cdots\otimes f_{m-1}^{}\rangle=
\begin{xy}
,(0,0)*\cir(10,10){},
,(0,0);(-7,-7)*{\bullet}**\crv{(-2,-4)} ?>*\dir{>}
,(0,0);(7,-7)*{\bullet}**\crv{(2,-4)} ?>*\dir{>}
,(0,0);(-4,-9)*{\bullet}**\crv{(-1,-4)} ?>*\dir{>}
,(-8,-9)*{\scriptstyle{f_0^{}}}
,(-4,-11)*{\scriptstyle{f^{}_1}}
,(1,-8)*{\cdots}
,(8,-9)*{\scriptstyle{f^{}_{m-1}}}
,(0,1.5)*{\scriptstyle{\xi}}
,(-4,-3)*{\scriptscriptstyle{i_0}}
,(-.5,-6)*{\scriptscriptstyle{i_1}}
,(5.6,-3)*{\scriptscriptstyle{i_{m-1}}}
\end{xy}\\
\\&=m!\,
\xi^{i_0i_1\dots i_{m-1}}(\partial_{i_0}f_0^{})
(\partial_{i_1}f_1^{})\cdots (\partial_{i_{m-1}}f_{m-1}^{}).
\end{align*}
It is an element of ${\mathcal G}^{0,1}$, independently of $m$.
\end{rem}

\subsection{A computation example}
 We want to compute the pairing
\[
\langle\U_{\Gamma}(\xi_1\otimes\xi_2)\vert f\otimes g\otimes
h\rangle,
\]
where $\Gamma$ is the graph of Figure 1 in \S \ref{SS:kag}.
The first step consists in decorating the boundary vertices with the functions $f,g$ and $h$, and the internal vertices with the polyvector fields $\xi_1$ and $\xi_2$:
\[
\begin{xy}
,(0,0)*\cir(10,10){},
,(-7,-7)*{\bullet},(7,-7)*{\bullet}
,(1,-3);(-7,-7)**\crv{(0,-4)} ?>*\dir{>}
,(1,-3);(0,-10)*{\bullet}**\crv{(0,-4)} ?>*\dir{>}
,(3,3);(1,-3)**\crv{(3,0)} ?>*\dir{>}
,(3,3);(7,-7)**\crv{(5,-5)} ?>*\dir{>}
,(3.8,4.6)*{\scriptstyle{\xi_1^{}}}
,(0,-1.4)*{\scriptstyle{\xi_2^{}}}
,(1.8,1.0)*{\scriptscriptstyle{1}}
,(-3,-4)*{\scriptscriptstyle{3}}
,(5.9,-2.5)*{\scriptscriptstyle{2}}
,(-1,-6.5)*{\scriptscriptstyle{4}}
,(-8,-9)*{\scriptstyle{f}}
,(0,-12)*{\scriptstyle{g}}
,(8,-9)*{\scriptstyle{h}}
\end{xy}
\]
Note that both the internal vertices of $\Gamma$ have
 two outgoing edges. This means that both $\xi_1$ and
$\xi_2$ have to be bi-vector fields in order to have a
non-trivial pairing. Let us write
$\xi_1=\xi_1^{i_1i_2}\partial_{i_1}\wedge\partial_{i_2}$ and
$\xi_2=\xi_2^{i_3i_4}\partial_{i_3}\wedge\partial_{i_4}$. We
now ``color'' the edges of $\Gamma$. The $k^{\rm th}$ edge is
labeled by $i_k$:
\[
\begin{xy}
,(0,0)*\cir(40,40){},
,(-8,-8)*{\bullet},(8,-8)*{\bullet}
,(1,-3);(-8,-8)**\crv{(0,-4)} ?>*\dir{>}
,(1,-3);(0,-11)*{\bullet}**\crv{(0,-4)} ?>*\dir{>}
,(3.5,4);(1,-3)**\crv{(3,0)} ?>*\dir{>}
,(3.5,4);(8,-8)**\crv{(5,-5)} ?>*\dir{>}
,(4.8,5.5)*{\scriptstyle{\xi_1}}
,(.1,-1)*{\scriptstyle{\xi_2}}
,(1.5,2)*{\scriptscriptstyle{i_1}}
,(-4.5,-4.5)*{\scriptscriptstyle{i_3}}
,(7,-2.5)*{\scriptscriptstyle{i_2}}
,(2.2,-6.9)*{\scriptscriptstyle{i_4}}
,(-9,-10)*{\scriptstyle{f}}
,(0,-13)*{\scriptstyle{g}}
,(9,-10)*{\scriptstyle{h}}
\end{xy}
\]
Now we are left with the task of giving a value to
the above colored graph. 
To do this we look at the values of the vertices:
\[
\begin{xy}
,(0,0);(-3,-3)**\dir{-}?>*\dir{>}
,(0,0);(3,-3)**\dir{-}?>*\dir{>}
,(0,1.5)*{\scriptstyle{\xi_1}}
,(-4,-2)*{\scriptscriptstyle{i_1}}
,(4,-2)*{\scriptscriptstyle{i_2}}
\end{xy}
\mapsto 2\,\xi_1^{i_1i_2}\,;\qquad\qquad
\begin{xy}
,(0,4);(0,0)**\dir{-} ?>*\dir{>}
,(0,0);(-3,-3)**\dir{-}?>*\dir{>}
,(0,0);(3,-3)**\dir{-}?>*\dir{>}
,(2,1)*{\scriptstyle{\xi_2}}
,(-4,-2)*{\scriptscriptstyle{i_3}}
,(4,-2)*{\scriptscriptstyle{i_4}}
,(-1,3.5)*{\scriptscriptstyle{i_1}}
\end{xy}
\mapsto2\, \partial_{i_1}\xi_2^{i_3i_4}\,;
\]
\begin{align*}
\begin{xy}
,(0,4);(0,-2)*{\bullet}**\dir{-} ?>*\dir{>}
,(0,-4)*{\scriptstyle{f}}
,(-1.5,3.5)*{\scriptscriptstyle{i_3}}
\end{xy}
\mapsto \partial_{i_3}f\,;\qquad\qquad
\begin{xy}
,(0,4);(0,-2)*{\bullet}**\dir{-} ?>*\dir{>}
,(0,-4)*{\scriptstyle{g}}
,(-1.5,3.5)*{\scriptscriptstyle{i_4}}
\end{xy}
\mapsto \partial_{i_4}g\,;\qquad\qquad
\begin{xy}
,(0,4);(0,-2)*{\bullet}**\dir{-} ?>*\dir{>}
,(0,-4)*{\scriptstyle{h}}
,(-1.5,3.5)*{\scriptscriptstyle{i_2}}
\end{xy}
\mapsto \partial_{i_2}h\,
\end{align*}
and multiply them. In conclusion, we get the formula
\[
\langle\U_{\Gamma}(\xi_1\otimes\xi_2)|f\otimes 
g\otimes h\rangle=4\,
\xi_1^{i_1i_2}(\partial_{i_1}\xi_2^{i_3i_4})(\partial_{i_3}f)
(\partial_{i_4}g)(\partial_{i_2}h).
\]

\section{The Chevalley-Eilenberg DGLA
$CE^\bullet(T_{\rm poly}^\bullet,D_{\rm
poly}^\bullet)$}\label{S:LA} We review the Lie algebra
structures of the main DGLAs involved, together with the
corresponding canonical differentials. This is standard material
included for the sake of completeness.

\subsection{The Lie algebra $T^\bullet_{\rm poly}$}
The space $T^\bullet_{\rm poly}$ of polyvector fields on
${\setR}^d$ is endowed with a graded Lie algebra
structure by the 
 Schouten-Nijenhuis bracket:
\begin{align*}
[\xi_0\wedge&\cdots\wedge\xi_k,\eta_0\wedge\cdots\wedge\eta_l]^{}_{SN}=\\
&=\sum_{i,j}(-1)^{i+j+k}[\xi_i,\eta_j]\wedge
\xi_0\wedge\cdots\wedge\hat{\xi}_i\wedge\cdots\wedge\xi_k\wedge
\eta_0\wedge\cdots\wedge\hat{\eta}_j
\wedge\cdots\wedge\eta_l,
\end{align*}
where $k,l\ge0$ and $\xi_i,\eta_j\in T_{\rm poly}^0$;
the hat over a factor denotes its absence. We will use the shorthand
notation
\[
\xi_{\bar i}:=\xi_{o}\wedge\cdots\wedge\hat\xi_i\wedge\cdots\wedge\xi_k
\]
so that the Schouten-Nijenhuis bracket
is written as
$$[\xi_0\wedge\cdots\wedge\xi_k,\eta_0\wedge\cdots\wedge\eta_l]^{}_{SN}=
\sum_{i=0}^k\sum_{j=0}^l(-1)^{i+j+k}[\xi_i,\eta_j]\wedge
\xi_{\bar i}\wedge{\eta}_{\bar j}.$$
The Schouten-Nijenhuis bracket is induced by a pre-Lie
operation on polyvector fields, the Nijenhuis-Richardson
pre-Lie multiplication  (\cite{K1}, p.15; \cite{Kel}, \cite{Mo} p.81). 
By this we mean that there is a (non associative)
multiplication $$
\bullet\colon T_{\rm poly}^k\otimes T_{\rm
poly}^l\to T_{\rm poly}^{k+l}
$$ such that (\cite{K1}, p.15):
$$[\xi,\eta]^{}_{SN}=\xi\bullet\eta-(-1)^{kl}
\eta\bullet\xi.$$
The definition of the Nijenhuis-Richardson
 $\bullet$ operation is the following:
if $\xi=\xi_0\wedge\cdots\wedge\xi_k$ and 
$\eta=\eta_0\wedge\cdots\wedge\eta_l$, with $\xi_i,\eta_j\in
T_{\rm poly}^0$, then
$$\xi\bullet\eta=\sum_{i,j}\pm\xi_i(\eta_j)
\wedge{\xi}_{\bar i}\wedge{\eta_{\bar j}},$$
where $\xi_i(\eta_j)=\xi_i^a(\partial_a\eta_j^b)\partial_b$. Note that,
if one writes $\xi=\xi^{i_0\dots
i_k}\partial_{i_0}\wedge\cdots\wedge\partial_{i_k}$ and
$\eta=\eta^{j_0\dots j_l}\partial_{j_0}\wedge\cdots\wedge\partial_{j_l}$,
 the formula of the Nijenhuis-Richardson
pre-Lie multiplication is written as (\cite{AMM}, p.19):
$$\xi\bullet\eta=\sum_{r=1}^{k}(-1)^{r-1}\xi^{i_0\ldots
i_k}(\partial_{i_r}\eta^{j_0\ldots j_l})
\partial_{i_0}\wedge\cdots\wedge\hat{\partial}_{i_r}\wedge\cdots\wedge
\partial_{i_k}\wedge
\partial_{j_0}\wedge\cdots\wedge\partial_{j_l},$$
i.e., the Nijenhuis-Richardson multiplication can be written as an alternate sum of
elementary multiplications $\bullet_r$, where
$$\xi\bullet_{r}\eta=\xi^{i_0\ldots
i_k}(\partial_{i_r}\eta^{j_0\ldots j_l})
\partial_{i_0}\wedge\cdots\wedge\hat{\partial}_{i_r}\wedge\cdots\wedge
\partial_{i_k}\wedge
\partial_{j_0}\wedge\cdots\wedge\partial_{j_l}.$$
 For instance, if $\alpha\in T^1_{\rm
poly}$ and
$\xi\in T^0_{\rm poly}$, then
\begin{align}
\alpha\bullet\xi&=\alpha^{i_1i_2}(\partial_{i_1}\xi^{i_3})
\partial_{i_2}\wedge\partial_{i_3}
-\alpha^{i_1i_2}(\partial_{i_2}\xi^{i_3})
\partial_{i_1}\wedge\partial_{i_3};\\
\xi\bullet\alpha&=\xi^{i_1}(\partial_{i_1}\alpha^{i_2i_3})\partial_{i_2}
\wedge\partial_{i_3}
\end{align}
Here is another example. If $\eta\in T^2_{\rm poly}$ and $h\in T^{-1}_{\rm
poly}$, then
\begin{align*}
\eta\bullet h&=3\eta^{i_1i_2i_3}(\partial_{i_1}h)\partial_{i_2}\wedge\partial_{i_3},
\\
h\bullet\eta&=0.
\end{align*}
The trivial differential $d=0$ makes $T^\bullet_{\rm poly}$ a differential graded
Lie algebra.
\begin{rem}
The above constructions can be globalized to a differential manifold $M$; in this
case the differential graded Lie algebra of polyvector fields is denoted by 
$$T_{\rm poly}^{\bullet}(M)=\bigoplus_{k=-1}^\infty T_{\rm
poly}^{k}(M)
=\bigoplus_{k=-1}^\infty H^0(M,\wedge^{k+1}TM).$$
\end{rem}

\subsection{The Lie algebra $D_{\rm poly}^\bullet$}\label{SS:Dpoly}

The Hochschild complex of the algebra ${\cal C}^\infty({\setR}^d)$ of smooth
functions on
${\setR}^d$ is a differential graded Lie algebra whose underlying graded vector
space is
\begin{align*}
{\rm Hoch}^\bullet({\cal C}^\infty({\setR}^d);{\cal
C}^\infty({\setR}^d))&=\bigoplus_{k=-1}^\infty {\rm Hoch}^k({\cal
C}^\infty({\setR}^d);{\cal C}^\infty({\setR}^d))\\
&=\bigoplus_{k=-1}^\infty{\rm
Hom}_{\setR}({\cal C}^\infty({\setR}^d)^{\otimes(k+1)},{\cal C}^\infty({\setR}^d)).
\end{align*}
The Lie bracket is the Gerstenhaber bracket: if $\varphi\in {\rm Hoch}^k({\cal
C}^\infty({\setR}^d);{\cal C}^\infty({\setR}^d))$ and $\psi\in {\rm Hoch}^l({\cal
C}^\infty({\setR}^d);{\cal C}^\infty({\setR}^d))$, then the bracket
$[\varphi,\psi]^{}_G$ is the element of ${\rm Hoch}^{k+l}({\cal
C}^\infty({\setR}^d);{\cal C}^\infty({\setR}^d))$ defined by
\begin{align*}
\langle &[\varphi,\psi]^{}_G|f_0^{}\otimes\cdots\otimes f_{k+l}^{}\rangle=
\\&=
\sum_{i=0}^k(-1)^{il}\langle\varphi|f_0^{}\otimes\cdots\otimes f_{i-1}^{}\otimes
\langle\psi|f_i^{}\otimes\cdots f_{i+l}^{}\rangle\otimes f_{i+l+1}\otimes\cdots
\otimes f^{}_{k+l}\rangle\\
&\phantom{m}-
\sum_{i=0}^l(-1)^{k(l+i)}\langle\psi|f_0^{}\otimes\cdots\otimes
f_{i-1}^{}\otimes
\langle\varphi|f_i^{}\otimes\cdots f_{i+k}^{}\rangle\otimes f_{i+k+1}\otimes\cdots
\otimes f^{}_{k+l}\rangle
\end{align*}
\par 
As the Schouten-Nijenhuis bracket, the Gerstenhaber bracket is also induced by a
pre-Lie multiplication, the Gerstenhaber composition
\begin{align*}
\circ\colon {\rm Hoch}^k({\cal C}^\infty({\setR}^d);{\cal
C}^\infty({\setR}^d))\otimes {\rm Hoch}^l({\cal
C}^\infty&({\setR}^d);{\cal
C}^\infty({\setR}^d))\to\\
&\to {\rm Hoch}^{k+l}({\cal
C}^\infty({\setR}^d);{\cal
C}^\infty({\setR}^d))
\end{align*}
defined by
\begin{align*}
\langle &\varphi\circ\psi|f_0^{}\otimes\cdots\otimes f_{k+l}^{}\rangle=
\\&=
\sum_{i=0}^k(-1)^{il}\langle\varphi|f_0^{}\otimes\cdots\otimes f_{i-1}^{}\otimes
\langle\psi|f_i^{}\otimes\cdots f_{i+l}^{}\rangle\otimes f_{i+l+1}\otimes\cdots
\otimes f^{}_{k+l}\rangle
\end{align*}
It is immediate to see that
\[
[\varphi,\psi]^{}_G=\varphi\circ\psi-(-1)^{kl}\psi\circ\varphi.
\]
Moreover, as for the Nijenhuis-Richardson multiplication $\bullet$, 
the Gerstenhaber multiplication $\circ$ can also be written as a graded alternate sum of
elementary compositions $\circ_i^{}$, where
\[
\langle\varphi\circ_i^{}\psi|f_0^{}\otimes\cdots\otimes f_{k+l}^{}\rangle=
\langle\varphi|f_0^{}\otimes\cdots\otimes f_{i-1}^{}\otimes
\langle\psi|f_i^{}\otimes\cdots f_{i+l}^{}\rangle\otimes f_{i+l+1}\otimes\cdots
\otimes f^{}_{k+l}\rangle
\]
 If we denote by $\mu$ the multiplication of functions in ${\cal C}^\infty(M)$, 
then $\mu$ is an element in ${\rm
Hoch}^1({\cal C}^\infty({\setR}^d);{\cal C}^\infty({\setR}^d))$ and the
associativity of $\mu$ is equivalent to $[\mu,\mu]^{}_G=0$. This means that $d_{\rm
Hoch}=[\mu,\cdot]^{}_G$ is a degree one differential of ${\rm Hoch}^\bullet({\cal
C}^\infty({\setR}^d);{\cal C}^\infty({\setR}^d))$. 
The differential $d_{\rm Hoch}$
will be called the Hochschild differential. Explicitly, if $\varphi\in
{\rm Hoch}^k({\cal C}^\infty({\setR}^d);{\cal
C}^\infty({\setR}^d))$, then $d_{\rm Hoch}\varphi$ is the element in 
${\rm Hoch}^{k+1}({\cal C}^\infty({\setR}^d);{\cal
C}^\infty({\setR}^d))$ defined by
\begin{align*}
\langle &d_{\rm Hoch}\varphi|f_0^{}\otimes\cdots\otimes f_{k+1}^{}\rangle=
\\&=
f_0^{}\langle\varphi\,|f_1^{}\otimes\cdots\otimes
f_{k+1}^{}\rangle-\sum_{i=0}^k\langle\varphi|f_0^{}\otimes\cdots
(f_{i}^{}\cdot f_{i+1}^{})\otimes\cdots\otimes f_{k+1}\rangle\\
&\phantom{mmm}+(-1)^k
\langle\varphi|f_0^{}\otimes\cdots\otimes
f_{k}^{}\rangle f^{}_{k+1}.
\end{align*}
 \par
 For additional details see
\cite{K1} \S 3.4.2 and \cite{CI} \S 3.2.2.
\par
One can define a cup-product
\begin{align*}
\cup\colon {\rm Hoch}^k({\cal C}^\infty({\setR}^d);{\cal
C}^\infty({\setR}^d))\otimes {\rm Hoch}^l({\cal
C}^\infty&({\setR}^d);{\cal
C}^\infty({\setR}^d))\to \\
&\to{\rm Hoch}^{k+l+1}({\cal
C}^\infty({\setR}^d);{\cal
C}^\infty({\setR}^d))
\end{align*}
as
\[
\langle D_1\cup D_2|f_0^{}\cdots f^{}_k\otimes f^{}_{k+1}\cdots
f_{k+l+1}^{}\rangle= \langle D_1|f_0^{}\otimes\cdots\otimes f^{}_k\rangle\cdot \langle
 D_2|f^{}_{k+1}\otimes\cdots\otimes f_{k+l+1}^{}\rangle
\]
The space of polydifferential operators on ${\setR}^d$ is the
smallest subspace $D^\bullet_{\rm poly}$ containing vector fields on
${\setR}^d$ and the multiplications by smooth functions on ${\setR}^d$, which is
closed under the cup product and the $\circ$ multiplication. 
More explicitly, elements of
$D^k_{\rm poly}$ are the linear operators on
${\cal C}^\infty({\setR}^d)^{\otimes (k+1)}$ with values in ${\cal
C}^\infty({\setR}^d)$ which can be written as
$$f_0\otimes\cdots\otimes f_k\mapsto C^{I_0,\dots,I_k}
(\partial_{I_0}f_0)\cdots(\partial_{I_k}f_k),$$
for suitable multi-indices $I_\bullet$ and smooth functions
$C^{I_\bullet}$. It is immediate to see that $D^\bullet_{\rm poly}$
is a differential graded Lie subalgebra of  ${\rm Hoch}^\bullet({\cal
C}^\infty({\setR}^d);{\cal C}^\infty({\setR}^d))$.  

\begin{rem}
More in general, given a smooth manifold $M$ one can consider the
differential
graded Lie algebra $D^\bullet_{\rm poly}(M)$ of polydifferential operators on $M$;
it is a sub differential graded Lie algebra of 
${\rm Hoch}^\bullet({\cal C}^\infty(M),{\cal C}^\infty(M))$.
\end{rem}

\subsection{The Lie algebra
$CE^\bullet(T^\bullet_{\rm poly},D^\bullet_{\rm poly})$}\label{SS:CE}

In this section we define the Chevalley-Eilenberg differential graded Lie algebra
$CE^\bullet(T_{\rm poly}^\bullet,D_{\rm poly}^\bullet)$. The
construction we are going to describe is just a particular case of
a more general construction of a differential graded Lie algebra
$CE^\bullet(\g_1,\g_2)$ out of two differential graded Lie algebras
$\g_1$ and $\g_2$, see \cite{Kel} for hints and details.
\par
As a bigraded vector space, 
\begin{align*}
CE^{\bullet,\bullet}(T^\bullet_{\rm poly},D^\bullet_{\rm
poly})&=\bigoplus_{p=-\infty}^\infty\bigoplus_{q=0}^\infty
CE^{p,q}(T^\bullet_{\rm poly},D^\bullet_{\rm poly})\\
&=
\bigoplus_{p=-\infty}^\infty\bigoplus_{q=0}^\infty{\rm
Hom}(\wedge^{q}T^\bullet_{\rm poly},D^\bullet_{\rm poly}[p]),
\end{align*}
where by $D^\bullet_{\rm
poly}[p]$ we mean the graded vector space $D^\bullet_{\rm
poly}$ with degrees shifted by $p$. The graded vector space 
$CE^\bullet(T^\bullet_{\rm poly},D^\bullet_{\rm
poly})$ is the graded vector space associated with the bigraded space
$CE^{\bullet,\bullet}$, i.e.,
\begin{align*}
CE^n(T^\bullet_{\rm poly},D^\bullet_{\rm
poly})&=\bigoplus_{p+q=n} CE^{p,q}(T^\bullet_{\rm poly},D^\bullet_{\rm
poly})\\
&=\bigoplus_{q=0}^\infty {\rm Hom}(\wedge^{q}T^\bullet_{\rm
poly},D^\bullet_{\rm poly}[n-q])
\end{align*}
The Lie
bracket
\begin{align*}
[\,,\,]^{}_{CE}\colon CE^{p_1,q_1}(T^\bullet_{\rm poly},D^\bullet_{\rm
poly})\otimes CE^{p_2,q_2}(&T^\bullet_{\rm poly},D^\bullet_{\rm
poly})\to\\
&\to CE^{p_1+p_2,q_1+q_2}(T^\bullet_{\rm poly},D^\bullet_{\rm
poly})
\end{align*}
is defined by
\begin{align*}
[{\mathcal
F},{\mathcal
H}]^{}_{CE}&(\gamma_1^{}\wedge\cdots\wedge\gamma_{q_1+q_2}^{})=\\
&\hskip-1em=\frac{1}{q_1!q_2!}\!\sum_{\sigma\in\Sigma_{q_1+q_2}}[{\mathcal
F}(
\gamma_{\sigma(1)}\wedge\cdots \wedge\gamma_{\sigma(q_1)}),
{\mathcal H}(
\gamma_{\sigma(q_1+1)}\wedge\cdots \wedge\gamma_{\sigma(q_1+q_2)})]^{}_G.
\end{align*}
Being bigraded, the Lie bracket defined above induces a graded Lie bracket on
$CE^\bullet(T^\bullet_{\rm poly},D^\bullet_{\rm
poly})$. If we look at the multiplication $\mu$ on smooth functions on ${\setR}^d$
as an element in $(D^\bullet_{\rm poly}[1])^0$, then the map ${\mathcal
F}_0:1\mapsto \mu$ is an element in $CE^{1,0}(T^\bullet_{\rm poly},D^\bullet_{\rm
poly})\subseteq CE^{1}(T^\bullet_{\rm poly},D^\bullet_{\rm
poly}) $. Since
\[
[{\mathcal F}^{}_0,{\mathcal F}^{}_0]_{CE}^{}(1)=[\mu,\mu]_{G}^{}=0,
\]
we have an {\em horizontal differential} 
\[
d_{1}^{}=[{\mathcal F}_0^{},\cdot]_{CE}^{}\colon
CE^{p,q}(T^\bullet_{\rm poly},D^\bullet_{\rm
poly})\to CE^{p+1,q}(T^\bullet_{\rm poly},D^\bullet_{\rm
poly})
\]
which makes $CE^{\bullet}(T^\bullet_{\rm poly},D^\bullet_{\rm
poly})$ a differential graded Lie algebra. Explicitly, if ${\mathcal F}\in CE^{p,q}
(T^\bullet_{\rm poly},D^\bullet_{\rm
poly}),$ then
\[
(d_1^{}{\mathcal F})(\gamma_1^{}\wedge\cdots\wedge\gamma_q^{})=d_{\rm
Hoch}({\mathcal F}(\gamma_1^{}\wedge\cdots\wedge\gamma_q^{})).
\]
We have not used the Lie algebra structure on $T^\bullet_{\rm poly}$, yet. By
looking at $D^\bullet_{\rm poly}$ as a trivial $T^\bullet_{\rm poly}$-module, we
can consider the classical Chevalley-Eilenberg differential as a vertical
differential on $CE^{\bullet,\bullet}(T^\bullet_{\rm poly},D^\bullet_{\rm poly})$.
Explicitly, the differential
\[
d_2^{}\colon CE^{p,q}(T^\bullet_{\rm poly},D^\bullet_{\rm
poly})\to CE^{p,q+1}(T^\bullet_{\rm poly},D^\bullet_{\rm
poly})
\]
is defined by
\[
(d_{2}^{}{\mathcal F})(\gamma_1\wedge\cdots\wedge\gamma_{q+1})=
\sum_{i<j}\pm{\mathcal F}([\gamma_i,\gamma_j]^{}_{SN}\wedge
\gamma_1\wedge\cdots
\wedge\widehat{\gamma_i}\wedge\cdots\wedge\widehat{\gamma_j}\wedge
\cdots\wedge\gamma_{q+1}^{}).
\]
By looking at an expression of the form
\[
\langle {\mathcal F}(\gamma_1\wedge\cdots\wedge \gamma_n)|f^{}_0\otimes\cdots
\otimes f_m^{}\rangle
\]
one sees that $d_1$ acts on the right side of the pairing whereas $d_2$ acts on the
left side. As a consequence, the two differentials commute, 
and define a degree one total differential 
$d_{CE}=d_1\pm d_2$ making $CE^\bullet(T^\bullet_{\rm
poly},D^\bullet_{\rm poly})$ a differential graded Lie algebra.

\section{The Maurer-Cartan equation and
$L_\infty$-morphisms}\label{S:mclinf}

Recall that Kontsevich's approach consists in pushing forward 
a solution of the Maurer-Cartan equation in
$T^\bullet_{\rm poly}$ (the Poisson structure $\alpha$) via
an $L_\infty$-morphism
${\mathcal F}:T^\bullet_{\rm poly}\to D^\bullet_{\rm poly}$,
 to get a solution of the Maurer-Cartan equation in
$D^\bullet_{\rm poly}$ (the star-product). In this section we
explain how
$L_\infty$-morphisms correspond  to solutions of the
Maurer-Cartan equation in the Chevalley-Eilenberg
differential graded Lie algebra $CE^\bullet(T^\bullet_{\rm
poly},D^\bullet_{\rm poly})$ and how each of these solutions
leads to an associative star-product quantizing a given
Poisson structure.
\par

\subsection{The Maurer-Cartan equation on
the Chevalley-Eilen\-berg complex} Given
a differential graded Lie algebra
$\g=\bigoplus_{i=0}^\infty\g^i$, the Maurer-Cartan
equation for $\g$ is the following equation for an element
$\gamma$ in $\g^1$:
\[
d\gamma+\frac{1}{2}[\gamma,\gamma]=0.
\]
The set of the solutions of the Maurer-Cartan equation for
$\g$ will be denoted by $MC(\g)$. We are interested in the set 
$MC(CE^\bullet(T^\bullet_{\rm poly},D^\bullet_{\rm poly}))$. 
Since 
\begin{align*}
CE^1(T^\bullet_{\rm poly},D^\bullet_{\rm
poly})&=\bigoplus_{p+q=1}CE^{p,q}(T^\bullet_{\rm
poly},D^\bullet_{\rm poly})\\
&=
\bigoplus_{q=0}^\infty{\rm Hom}(\wedge^qT^\bullet_{\rm
poly},D^\bullet_{\rm poly}[1-q]),
\end{align*}
a solution of the Maurer-Cartan equation for the
Chevalley-Eilenberg complex of $T^\bullet_{\rm
poly}$ and $D^\bullet_{\rm
poly}$ is 
\[
{\mathcal F}={\mathcal F}_0^{}+{\mathcal
F}_1^{}+{\mathcal F}_2^{}+\cdots,
\]
where
\[
{\mathcal F}_q^{}\colon \wedge^qT^\bullet_{\rm poly}\to
D^\bullet_{\rm poly}[1-q]
\]
and the ${\mathcal F}_q^{}$'s satisfy the equations
\[
d_1^{}{\mathcal
F}_{q}^{} \pm d^{}_2{\mathcal
F}_{q-1}^{}+\frac{1}{2}\sum_{q_1+q_2=q}
[{\mathcal F}_{q_1}^{},{\mathcal F}_{q_2}^{}]_{CE}^{}=0, \qquad q\ge 0,
\]
corresponding to the ${\rm Hom}(\wedge^q T^\bullet_{\rm
poly},D^\bullet_{\rm poly}[2-q])$-components of the single
Maurer-Cartan equation
\[
d_{CE}^{}{\mathcal
F}+\frac{1}{2}
[{\mathcal F},{\mathcal F}]_{CE}^{}=0.
\]
 In particular
${\mathcal F}_0^{}$ has to satisfy the equation
$d_{\rm Hoch}({\mathcal F}_0^{}(1))+\frac{1}{2}[{\mathcal
F}_0^{}(1),{\mathcal F}_0^{}(1)]_{G}^{}=0$, that is $\mu+{\mathcal
F}_0^{}(1)$ has to be an associative product on the space of smooth
functions on ${\setR}^d$. This equation is trivially satisfied if
we choose
${\mathcal F}_0^{}=0$: as usual, one should think of
solutions of the Maurer-Cartan equations as {\em
perturbations}; when added to the {\em base point} (the
``initial value'') they yield the desired {\em deformations}.

If we
consider the space
$MC_{0}(CE^\bullet(T^\bullet_{\rm poly},D^\bullet_{\rm
poly}))$ of solutions of the Maurer-Cartan equation on the
Chevalley-Eilenberg complex of $T^\bullet_{\rm poly}$ and
$D^\bullet_{\rm poly}$ such that the ``initial
value'' ${\mathcal F}_0^{}$ is zero, then the equations for
an element ${\mathcal F}={\mathcal F}_0^{}+{\mathcal
F}_1^{}+{\mathcal F}_2^{}+\cdots$ of 
$MC_0(CE^\bullet(T^\bullet_{\rm poly},D^\bullet_{\rm poly}))$
are 
\begin{align*}
\begin{cases}
{\mathcal F}_0^{}=0\\
\\
\displaystyle{d_1^{}{\mathcal F}_q^{}\pm d^{}_2{\mathcal
F}_{q-1}^{}+\frac{1}{2}\sum_{\substack{q_1+q_2=q\\q_1,q_2
\geq1}}
[{\mathcal F}_{q_1}^{},{\mathcal F}_{q_2}^{}]_{CE}^{}=0,
\quad q\geq 1}.
\end{cases}
\end{align*}
Written out explicitly for $q$ polyvector fields
$\gamma^{}_1,\dots \gamma^{}_q$, these equations are
\begin{multline*}
d^{}_{\rm Hoch}({\mathcal
F}_{q}(\gamma_1^{}\wedge\cdots\wedge
\gamma_q^{}))+\sum_{i<j}\pm{\mathcal
F}_{q-1}([\gamma_i,\gamma_j]_{SN}\wedge
\gamma_1\wedge\cdots\wedge
\widehat{\gamma_i}\wedge\cdots\wedge\widehat{\gamma_j}\wedge
\cdots\wedge\gamma_q)+\\
+\frac{1}{2}\sum_{\substack{q_1+q_2=q\\q_1,q_2\geq1}}
\frac{1}{q_1!q_2!}\sum_{\sigma\in\Sigma_q}\pm[{\mathcal
F}_{q_1}(
\gamma_{\sigma(1)}\wedge\cdots \wedge\gamma_{\sigma(q_1)}),
{\mathcal F}_{q_2}(
\gamma_{\sigma(q_1+1)}\wedge\cdots\wedge
\gamma_{\sigma(q)})]_G^{}=0,
\end{multline*}
where $\Sigma_q$ denotes the symmetric group on $q$ elements.
But this precisely means that the maps ${\mathcal F}_q\colon\wedge^q
T^\bullet_{\rm poly}\to D^\bullet_{\rm poly}[1-q]$, $q\geq1$
are the Taylor coefficients of an $L_\infty$-morphism between the
differential graded
Lie algebras $T^\bullet_{\rm poly}$ and $D^\bullet_{\rm
poly}$. We have therefore proved the following
equivalence.
\begin{prop} The set ${\rm
Hom}_{L_\infty}(T^\bullet_{\rm poly},D^\bullet_{\rm poly})$
is naturally identified with the set
$MC_0^{}(CE^\bullet(T^\bullet_{\rm poly},D^\bullet_{\rm
poly}))$  of solutions of the Maurer-Cartan equation in
the Chevalley-Eilenberg differential graded Lie algebra
$CE^\bullet(T^\bullet_{\rm poly},D^\bullet_{\rm
poly})$ of $T^\bullet_{\rm poly}$ and $D^\bullet_{\rm
poly}$, such that
${\mathcal F}_0^{}=0$.
\end{prop}

\subsection{The quantization of Poisson structures}
Let now ${\mathcal F}$ be an element of $MC_0^{}(CE^\bullet(T^\bullet_{\rm
poly},D^\bullet_{\rm
poly}))$, and let
$\alpha$ be a Poisson structure on $\setR^d$. Then
$\alpha$ is an element of $T^1_{\rm poly}$ and
$[\alpha,\alpha]^{}_{SN}=0$. Since $\alpha$ lies in
$T^1_{\rm poly}$, ${\mathcal F}_q^{}(\alpha^{\wedge q})$ is
an element of $D^1_{\rm poly}$ for any $q\geq 0$. It follows
that the formal power series
\[
\pi_\alpha^{}(\hbar):={\mathcal F}(\exp\{\hbar \alpha\})
\]
is an element of $D^1_{\rm poly}[[\hbar]]$. Choosing
$\gamma_i=\alpha$ for any $i$ from 1 to $q$ in the equations
for the
${\mathcal F}^{}_q$'s written above, we find
\[
d^{}_{\rm Hoch}({\mathcal
F}_{q}(\alpha^{\wedge q}/{q!}))
+\frac{1}{2}\sum_{\substack{q_1+q_2=q\\q_1,q_2\geq1}}
\pm[{\mathcal
F}_{q_1}(\alpha^{\wedge q_1}/{q_1!}),
{\mathcal F}_{q_2}(\alpha^{\wedge q_2}/{q_2!})]_G^{}=0.
\]
Therefore
\[
d^{}_{\rm Hoch}\pi_\alpha^{}+\frac{1}{2}[\pi_\alpha^{},
\pi_\alpha^{}]_G^{}=0.
\]
In other words,
\[
\mu_\alpha^{}(\hbar):=\mu+\pi_{\alpha}^{}
\]
is an associative product (depending on the formal parameter
$\hbar$) on the space of the smooth functions on $\setR^d$.
Note that $\lim_{\hbar\to 0}\mu_\alpha^{}(\hbar)=\mu$, i.e.,
$\mu_\alpha^{}(\hbar)$ is an associative deformation of
$\mu$. For two smooth functions $f$ and $g$ on $\setR^d$, we
set
\[
f\star_\alpha^{}g:=\langle \mu_\alpha^{}(\hbar)|f\otimes
g\rangle=fg+\hbar\langle {\mathcal F}_1^{}(\alpha)|
f\otimes
g\rangle+O(\hbar^2).
\]
As an immediate consequence, we obtain
\begin{prop}
Any solution
${\mathcal F}$ of the Maurer-Cartan equation on the
Che\-val\-ley-Eilenberg
DGLA
$CE^\bullet(T^\bullet_{\rm poly},\allowbreak D^\bullet_{\rm
poly})$ with
${\mathcal F}_0^{}=0$ and
${\mathcal F}_1^{}(\gamma)=\gamma$ (with the usual
identification of
polyvector fields with polydifferential operators) leads to
a
\emph{star product} $\star_\alpha$ such that
\[
f\star_\alpha^{}g=fg+\hbar\{f,g\}_\alpha^{}
+O(\hbar^2),
\]
i.e., to a deformation quantization of the Poisson structure
$\alpha$.
\end{prop}

\subsection{The vertical Maurer-Cartan equation}
It is convenient to consider the translation of ${\mathcal
F}$ by $\mu$, since the defining equations for
${\mathcal F}^\mu_{}=\mu+{\mathcal F}$ are simpler than the
original equations for ${\mathcal F}$. In particular,
${\mathcal F}^\mu_{}$ can be seen as a solution of the
``vertical'' Maurer-Cartan equation in the
Chevalley-Eilenberg complex, i.e., ${\mathcal F}^\mu_{}$
satisfies
\[
d_2^{}{\mathcal F}^\mu_{}+\frac{1}{2}[{\mathcal
F}^\mu_{},{\mathcal F}^\mu_{}]_{CE}^{}=0
\]
Indeed, the components of ${\mathcal F}^\mu_{}$ are given by
\[
\begin{cases}
{\mathcal F}^\mu_{0}(1)=\mu+{\mathcal F}^{}_{0}(1);\\
{\mathcal F}^\mu_{q}={\mathcal F}_{q}^{}; \text{ for }q\geq 1,
\end{cases}
\]
 and we can rewrite the equations for
${\mathcal F}$ as
\begin{align*}
\begin{cases}
{\mathcal F}_0^{\mu}(1)=\mu;\\
\\
\displaystyle{d_1^{}{\mathcal F}_q^{\mu}\pm d^{}_2{\mathcal
F}_{q-1}^{\mu}+\frac{1}{2}\sum_{\substack{q_1+q_2=q\\q_1,q_2
\geq1}}
[{\mathcal F}_{q_1}^{\mu},{\mathcal
F}_{q_2}^{\mu}]_{CE}^{}=0,
\quad q\geq 1},
\end{cases}
\end{align*}
where we have used the fact that $d_2{\mathcal F}^\mu_0$ is
trivially zero since the ``vertical''
differential $d_2^{}$ is zero on $CE^{p,0}(T^\bullet_{\rm
poly},D^\bullet_{\rm poly})$. Since the ``horizontal''
differential
$d_1$ on the Chevalley-Eilenberg complex coincides with the
Lie-adjoint of
${\mathcal F}^\mu_{0}$, we can rewrite the above equations as
\begin{align*}
\begin{cases}
{\mathcal F}_0^{\mu}(1)=\mu;\\
\\
\displaystyle{ d^{}_2{\mathcal
F}_{q-1}^{\mu}+\frac{1}{2}\sum_{\substack{q_1+q_2=q\\q_1,q_2
\geq0}}
[{\mathcal F}_{q_1}^{\mu},{\mathcal
F}_{q_2}^{\mu}]_{CE}^{}=0,
\quad q\geq 1},
\end{cases}
\end{align*}
which is precisely the vertical Maurer-Cartan equation on 
$CE^{\bullet}(T^\bullet_{\rm
poly},D^\bullet_{\rm poly})$. 
\begin{rem} \label{R:defq}
The space 
$MC_\mu(CE^\bullet(T^\bullet_{\rm poly},
D^\bullet_{\rm poly}),d_2^{})$ 
of solutions of the vertical Maurer-Cartan
equation with ${\mathcal F}^\mu_0(1)=\mu$ should be thought
as an affine version of 
$MC_0(CE^\bullet(T^\bullet_{\rm poly},
D^\bullet_{\rm poly}),d_{CE}^{})$. 
In particular, the translation by $\mu$ establishes
a one to one correspondence between 
$MC_0(CE^\bullet(T^\bullet_{\rm poly},
D^\bullet_{\rm poly}),d_{CE}^{})$ 
and $MC_\mu(CE^\bullet(T^\bullet_{\rm poly},
D^\bullet_{\rm poly}),\allowbreak d_2^{})$. 
As a consequence, we get a
star-product quantization of Poisson structures out of any
element of 
$MC_\mu(CE^\bullet(T^\bullet_{\rm poly},D^\bullet_{\rm poly}),d_2^{})$ 
such that ${\mathcal F}^\mu_1\colon T^m_{\rm poly}\to D^m_{\rm poly}$ 
is the canonical identification of polyvector fields
with polydifferential operators.
\end{rem}

\section{The differential graded Lie algebra of graphs}\label{S:graphs}
Kontsevich's proof from \cite{K1} clearly hints to a coproduct structure 
controlling the coefficients of the formality morphism \cite{pqft}.
The corresponding result was stated in \cite{cfg}
in a framework dual to the DGLA picture briefly mentioned in \cite{comb}.
Here we focus on the DGLA structure with its two pre-Lie operations inspired 
by Connes-Kreimer's insertion operation and Gerstenhaber composition.

Recall from section \ref{SS:kag} that 
${\mathcal G}^{p,q}$ is the set of isomorphism classes of all
admissible graphs $\Gamma$ with
$q$ internal vertices and such that $v(\Gamma)-e(\Gamma)-1=p$, where
$v(\Gamma)$ and $e(\Gamma)$ are the number of vertices and edges of $\Gamma$,
respectively. Finally,
$${\mathscr G}^{\,\bullet,\bullet}
=\bigoplus_{\,p,q=0}^\infty{\mathscr G}^{p,q}$$ 
will denote
the bigraded vector space generated by
${\mathcal G}^{\bullet,\bullet}=\cup_{p,q=0}^\infty{\mathcal G}^{p,q}$.

If $\Phi$ is a
subgraph of $\Gamma$, we can collapse $\Phi$
to a point, obtaining a new graph $\Gamma/\Phi$; the point
corresponding to the collapsed subgraph will be a vertex
$v_{\Phi}^{}$ in the quotient graph, which can be 
either an internal or a boundary vertex. To describe
more precisely the collapsing of subgraphs, recall that
admissible graphs have to be thought as embedded into a
disk. To collapse a subgraph
$\Phi$ to a point, one first chooses a point $v_\Phi^{}$ in
the disk; then one considers a simply connected region
$R_\Phi$ in the disk having a connected (possibly
empty) intersection with the boundary of the disk and which
contains $v_\Phi$ and all the vertices of
$\Phi$, but contains no vertex of $\Gamma$ which is not in
$\Phi$. Note that this requirement implies that not every
subgraph of $\Gamma$ can be contracted to a point: for
instance, if $\Phi$ contains two boundary vertices $v_1$
and $v_2$ then, to be a contractible subgraph, $\Phi$ has to
contain also all the boundary vertices of $\Gamma$ which lie
between $v_1$ and $v_2$. Then one collapses the whole region
$R_\Phi$ to the point $v_\Phi^{}$
is such a way that, during the collapsing, the points of the
region which were on the boundary of the disk 
remain constrained on the boundary of the disk.  
Note that this implies that
$v_\Phi$ is forced to be a boundary vertex if $\Phi$ contains
at least a boundary vertex. 
\par
The set of admissible graphs is
not closed under the operation of quotienting by
contractible subgraphs, i.e., given an admissible graph
$\Gamma$ and a contractible subgraph
$\Phi$ of
$\Gamma$, the quotient graph
$\Gamma/\Phi$ is not necessarily an element of ${\mathcal
G}^{\bullet,\bullet}$. For instance, in the graph
$\Gamma/\Phi$ there could be a loop or an edge stemming from
a boundary vertex. 
\begin{definition}
We say
that the subgraph $\Phi$ is a \emph{normal subgraph} of $\Gamma$ if
$\Gamma/\Phi$ is an element of ${\mathcal G}^{\bullet,\bullet}$; 
in such a case we write $\Phi\normal\Gamma$. 
\end{definition}

\begin{rem}
 Note that a normal subgraph $\Phi\normal\Gamma$ 
has to be a \emph{full
subgraph}, i.e. every edge of $\Gamma$
connecting two vertices of $\Phi$ is an edge of $\Phi$. Indeed, if $\Phi$ were not full, then there
would be two vertices
$v_1$ and $v_2$ in $\Phi$ and an edge $e$ joining $v_1$ and $v_2$ in $\Gamma$
which is not an edge in $\Phi$. Then, in the quotient graph
$\Gamma/\Phi$ the edge $e$ would be a loop starting and
ending on the vertex $v_\Phi^{}$, and so
$\Gamma/\Phi$ would not be admissible.
\end{rem}
\begin{definition}
 If $\Phi$ is a normal subgraph of $\Gamma$ and $\Gamma/\Phi=\Psi$, we say
that
$\Gamma$ is an
\emph{extension} of $\Psi$ by $\Phi$ and write
\[
\Phi\hookrightarrow\Gamma\twoheadrightarrow\Psi
\]
\end{definition}
This simply means that $\Gamma$ is obtained by inserting the graph $\Phi$
into a vertex of $\Psi$. Depending on the type of vertex $\Phi$ is inserted
into, the extension is called internal or boundary extension.
\par
Note that given two graphs $\Phi$ and $\Psi$ we cannot talk of ``the
extension'' of $\Psi$ by $\Phi$ since there are many possible extensions,
even if one fixes the vertex of $\Psi$ where to insert $\Phi$. On the other
hand, the sum of all possible extensions of $\Psi$ by $\Phi$ is a
well-defined element of the vector space ${\mathscr G}^{\,\bullet,\bullet}$, so we
get a bilinear operation on ${\mathscr G}^{\,\bullet,\bullet}$. Depending whether
we consider all extensions or only internal/boundary extensions, we obtain the
Connes-Kreimer product
\[
\Psi *
\Phi=\sum_{\Phi\hookrightarrow\Gamma\twoheadrightarrow\Psi}\pm
\Gamma
\]
or the Hochschild-Kontsevich products
\[
\Psi \bullet
\Phi=\sum_{\substack{\Phi\hookrightarrow\Gamma\twoheadrightarrow\Psi\\
\text{internal}}}\pm
\Gamma
\]
\[
\Psi \circ
\Phi=\sum_{\substack{\Phi\hookrightarrow\Gamma\twoheadrightarrow\Psi\\
\text{boundary}}}\pm
\Gamma
\]
The signs in the above formulas depend upon the combinatorial
data of the graphs involved. Explicit formulas can be found,
for instance, in \cite{AMM}.
Note that the product $\circ$ above is bigraded:
\[
\circ\colon {\mathscr G}^{\,p_1,q_1}\otimes {\mathscr
G}^{\,p_2,q_2}\to {\mathscr G}^{\,p_1+p_2,q_1+q_2},
\]
whereas the $\bullet$ product has degree $(0,-1)$:
\[
\bullet\colon {\mathscr G}^{\,p_1,q_1}\otimes {\mathscr
G}^{\,p_2,q_2}\to {\mathscr G}^{\,p_1+p_2,q_1+q_2-1}.
\]
The Connes-Kreimer product $*$, being a combination of $\circ$ and $\bullet$, is
non homogeneous. In the rest of this
paper we will be interested only in the Hochschild-Kontsevich products
and the Lie algebra structures they induce; 
the reader interested in Connes-Kreimer product is referred to \cite{CK}.
\par

\subsection{The $\circ$ composition of graphs}
It is a simple check to see that the $\circ$ operation makes
${\mathscr G}^{\,\bullet,\bullet}$ a pre-Lie algebra, i.e. the (graded)
commutator defines a bigraded Lie algebra structure on
${\mathscr G}^{\,\bullet,\bullet}$.
\par
As for the Gerstenhaber composition, also the graph composition $\circ$ can be
written as an alternate sum of elementary compositions 
\begin{align*}
\circ_i^{}\colon {\mathscr G}^{\,p_1,q_1}\otimes {\mathscr
G}^{\,p_1,q_1}&\to {\mathscr G}^{\,p_1+p_2,q_1+q_2}\\
\Psi\otimes\Phi&\mapsto\Psi\circ_i^{}\Phi,
\end{align*}
which insert the graph $\Phi$ in the
$i^{\rm th}$ boundary vertex of the graph $\Psi$ (``a la Kathotia''
\cite{Kath}, p.15: a leg of $\Psi$ may land on any vertex of
$\Phi$; as we will see in the next section this corresponds to
Leibniz rule for  derivations).
As far as concerns the labeling of the edges and vertices of 
$\Psi\circ\Phi$,
the labeling of internal vertices of
$\Psi\circ_i\Phi$ is defined as the concatenation of the
corresponding labelings, the  labeling of boundary vertices of
$\Psi\circ_i\Phi$ is defined ``inserting'' the boundary labels of
$\Phi$ shifted by $i-1$ in place of the $i^{\rm th}$ boundary
vertex of
$\Psi$ and consequently shifting the labels on the other boundary
vertices of
$\Psi$; the total order on edges from each internal vertex is
unchanged.  Here is an example.

\vbox{\hskip -1.5 cm
\includegraphics{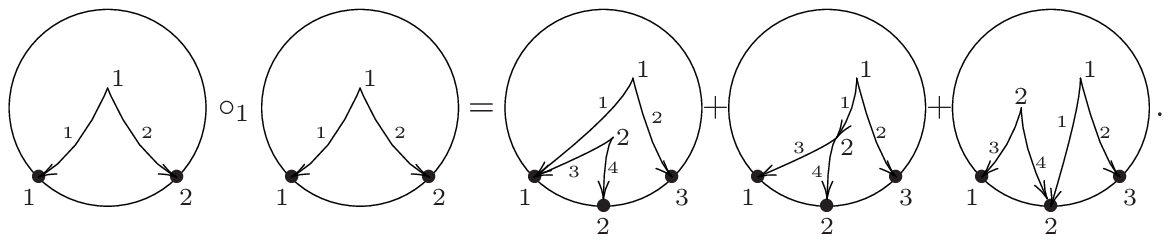}
}
 We will denote by $[\,,\,]$ the corresponding Lie bracket
\[
[\Psi,\Phi]=\Psi\circ\Phi\pm\Phi\circ\Psi
\]
and will call it the Hochschild-Kontsevich Lie bracket. 
Note that the graph $b_{0,2}$ with two boundary vertices 
and no internal vertices satisfies $[b_{0,2},b_{0,2}]=0$. 
Therefore it induces a differential $d_1=[b_{0,2},\cdot]$ on
${\mathscr G}^{\,\bullet,\bullet}$. 
Since $b_{0,2}\in{\mathscr G}^{\,1,0}$ 
and the Hochschild-Kontsevich Lie bracket is bigraded, $d_1^{}$ is a
degree $(1,0)$ differential:
\[
d_1^{}\colon {\mathscr G}^{\,p,q}\to {\mathscr G}^{\,p+1,q}.
\]

\subsection{The differential $d_2^{}$ on  graphs}
The ``vertical'' differential $d_2$ on admissible graphs is the usual 
graph homology differential introduced by Kontsevich in \cite{Kfc}, p.151
(see also \cite{Vor}, p.3): $d_2^{}\Psi$ is the alternate sum over all
graphs $\Gamma$ which can be obtained by expanding the internal vertices
of $\Psi$ by the insertion of an additional edge. In formulas,
 \[
d_2^{}\Psi=\sum_{\substack{e\hookrightarrow\Gamma\twoheadrightarrow\Psi\\
\text{internal}}}\pm
\Gamma
\]
This is conveniently rewritten in terms of the $\bullet$ multiplication
we described above:
\[
d_2^{}\Psi=\Psi\bullet e
\]
The usual argument shows that $d_2^{}$ is indeed a differential:
$(\Psi\bullet e)\bullet e$ is basically the sum over all ways of inserting two
edges in the internal vertices of $\Psi$, and each term in this sum
occurs twice with opposite signs. Note that, since the
$\bullet$ multiplication has degree $(0,-1)$ and $e\in{\mathscr
G}^{\,0,2}$, the differential $d_2^{}$ is a degree $(0,1)$ differential. 
\par As the $\circ$ multiplication, also the $\bullet$ multiplication can be
written as an alternate sum of elementary $\bullet_i^{}$ operations, which
insert
the second graph in the $i^{\rm th}$ internal vertex of the first graph. This
remark leads to the description of the differential $d_2{}$ acting on admissible
graphs as
$d^{}_2=\sum_i\pm d_{2,i}^{}$, where
$d_{2,i}^{}\colon{\mathscr G}^{\,p,q}\to{\mathscr G}^{\,p,q+1}$ is the operator
which inserts an edge (both ways) splitting the
$i^{\rm th}$ vertex of a graph. More precisely, $d_{2,i}\Psi=\sum\pm\Gamma$, where
the sum ranges over the set of graphs $\Gamma$ which are obtained by inserting an
edge (both ways) in the $i^{\rm th}$ vertex of
$\Psi$. 
\par
Here is an example.

\vbox{\hskip -1.7 cm
\includegraphics{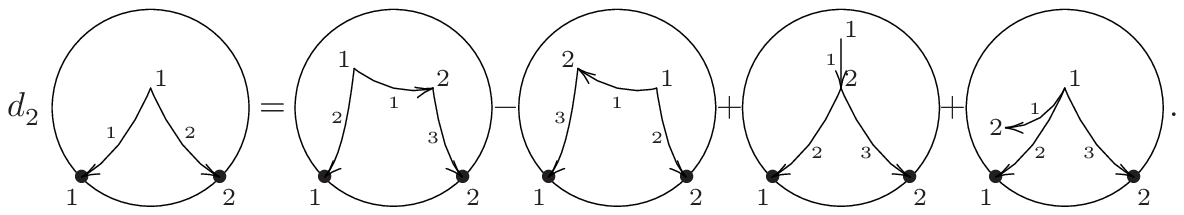}
}

\subsection{The total differential on graphs}
We have shown that there are two differentials on the bi-graded vector space
${\mathscr G}^{\,\bullet,\bullet}$. Namely, we have the degree $(1,0)$ differential
$d_1$ an the degree $(0,1)$ differential $d_2$. We will now show that the two
differentials commute, so that they induce a degree one total differential $d=d_1\pm
d_2$ on the total graph complex ${\mathscr G}^{\,\bullet}$, where
\[
{\mathscr G}^{\,n}=\bigoplus_{p+q=n}{\mathscr G}^{\,p,q}
\] 
Moreover, the total differential $d$ will turn out to be compatible with the graded
Lie bracket $[\,,\,]$ on ${\mathscr G}^{\,\bullet}$, so that we obtain a
differential graded Lie algebra of admissible graphs. We summarize these results in
the following
\begin{prop}
The Lie bracket $[\,,\,]$ and the differential $d$ endow the graded vector space
${\mathscr G}^{\,\bullet}$ with a differential graded Lie algebra structure, which we will
call the Hochschild-Kontsevich differential graded Lie algebra of graphs.
\end{prop}
\begin{proof}
Since the differential $d_1^{}$ is defined by means of the Lie bracket $[\,,\,]$ on
${\mathscr G}^{\,\bullet}$, and $d_2\,b_{0,2}=0$, the compatibility
of the differential
$d_2^{}$ with the Lie bracket will imply at once both the
compatibility of the two differentials and the compatibility of
the total differential with the Lie bracket.  Since the internal
vertices of
$\Psi\circ\Phi$ are exactly the disjoint union of the internal
vertices of $\Phi$ and of the internal vertices of
$\Psi$, we have 
\[
(\Psi\circ\Phi)\bullet e=(\Psi\bullet
e)\circ\Phi+\Psi\circ(\Phi\bullet e),
\]
which immediately implies the compatibility of $d_2{}$ with the Lie bracket.
\end{proof}
 
\begin{rem} The DGLA structure on graphs we have described is
actually part of a richer algebraic structure. Indeed, juxtaposition
of graphs defines a degree $(-1,0)$ operation $\cup\colon{\mathscr
G}^{\,p_1,q_1}\otimes {\mathscr G}^{\,p_2,q_2}\to{\mathscr
G}^{\,p_1+p_2-1,q_1+q_2}$ which, together with the Lie bracket
$[\,,\,]$, satisfies the axioms of Gerstenhaber algebra, up to
homotopy. On the
$T^\bullet_{\rm poly}$-side, the $\cup$ multiplication corresponds to
the usual wedge product of polyvector fields; on the $D^\bullet_{\rm
poly}$-side it corresponds to the cup-product on the Hochschild
complex. See \cite{gamella-halbout} and \cite{ginot-halbout} for
details.
\end{rem}

\section{The map $\U$ as a morphism of differential gra\-ded Lie
algebras}\label{S:Umap}

Recall from Section \ref{S:gc} that the graphical calculus for 
derivations can
be seen as a map
\[
\U\colon {\mathcal G}^{p,q}\to{\rm Hom}\left((T^\bullet_{\rm
poly})^{\otimes q},D_{\rm
poly}^\bullet[p]\right),
\]
Therefore, if we identify the antisymmetric product $\bigwedge^q T_{\rm
poly}^\bullet$ with a subspace of the tensor
product $\bigotimes^q T_{\rm poly}^\bullet$ via the usual
anti-symmetrization map, we can look at $\U$ as a map of bigraded vector spaces
\[
\U\colon {\mathscr G}^{\,p,q}\to CE^{p,q}(T^\bullet_{\rm
poly},D_{\rm
poly}^\bullet).
\]
We will show in this section that $\U$ is actually a map of differential graded Lie
algebras.

\begin{thm}\label{T:DGLA}
The graphical calculus for derivations $\U$ 
is a morphism of differential graded Lie algebras.
\end{thm}
The proof will be divided into two parts.

\subsection{$\U$ as a morphism of Lie algebras}
\begin{lem}
$[\U_{\Gamma_1},\U_{\Gamma_2}]_{CE}^{}=\U_{[\Gamma_1,\Gamma_2]}$
\end{lem}
\begin{proof} Let $\Gamma_1\in{\mathscr G}^{\,p_1,q_1}$ and $\Gamma_2\in{\mathscr
G}^{\,p_2,q_2}$. The lemma will immediately follow from the finer statement
$$\C{U}_{\Gamma_1}(\xi_1\otimes\cdots\otimes\xi_{q_1})\circ_i 
\C{U}_{\Gamma_2}(\xi_{q_1+1}\otimes\cdots\otimes\xi_{q_1+q_2})=
\C{U}_{\Gamma_1\circ_i\Gamma_2}(\xi_1\otimes\cdots\otimes\xi_{q_1+q_2}),$$
where the $\circ_i$ operation on the left is the $i^{\rm th}$ insertion in the
Gerstenhaber pre-Lie operation $\circ$, and the $\circ_i$ operation the right is
the operation that inserts the graph $\Gamma_2$ in the $i^{\rm th}$ boundary vertex
of the graph $\Gamma_1$. 
By definition of the
$\circ_i$ Gerstenhaber composition, what we have to show is that
\begin{align*}
\langle
\U^{}_{\Gamma_1\circ_i\Gamma_2}&(\xi_1\otimes\cdots\otimes
\xi_{q_1+q_2})|f_0\otimes\cdots\otimes f_{m_1+m_2-2}\rangle\\
&=
\langle\U_{\Gamma_1}(\xi_1\otimes\cdots\otimes 
\xi_{q_1})|f_0\otimes\cdots\otimes
f_{i-1}\otimes\\
&\phantom{mm}\otimes\langle \U_{\Gamma_2}(\xi_{q_1+1}\otimes\cdots\otimes
\xi_{q_1+q_2})| f_i\otimes\cdots\otimes f_{i+m_2-1}\rangle\otimes\\
 &\phantom{mm}\otimes f_{i+m_1}\otimes
\cdots\otimes f_{m_1+m_2-2}\rangle,
\end{align*}
where $m_1$ and $m_2$ are the number of boundary vertices of $\Gamma_1$ and
$\Gamma_2$, respectively. This directly follows from  the
definition of the $\circ_i$ graph boundary insertion and by
the Leibniz rule for derivations.
\end{proof}
\begin{example}
We will check that $b_{1,2}\circ_1 b_{1,2}$ is mapped by $\U$ to
$\U_{b_{1,2}}\circ_1 \U_{b_{1,2}}$, where $b_{1,2}$ is the 2-corolla
\[
\begin{xy}
,(0,0)*\cir(10,10){},
,(-7,-7)*{\bullet},(7,-7)*{\bullet}
,(0,2);(-7,-7)**\crv{(-3,-5)} ?>*\dir{>}
,(0,2);(7,-7)**\crv{(3,-5)} ?>*\dir{>}
,(1,3)*{\scriptstyle{1}}
,(-4,-2.5)*{\scriptscriptstyle{1}}
,(4,-2.5)*{\scriptscriptstyle{2}}
\end{xy}
\] 
Clearly, $\U_{b_{1,2}}(\xi_1), \U_{b_{1,2}}(\xi_2)$ and
$\U_{b_{1,2}\circ_1 b_{1,2}}(\xi_1\otimes \xi_2)$ vanish unless
$\deg\xi_1=\deg\xi_2=1$. 
So we are left to check that
$\U_{b_{1,2}}(\xi_1)\circ_1\U_{b_{1,2}}(\xi_2)=\U_{b_{1,2}\circ_1
b_{1,2}}(\xi_1\otimes
\xi_2)$ when both $\xi_1$ and $\xi_2$ are 2-vector fields. We begin by computing
$\U_{b_{1,2}}(\xi_1)\circ_1\U_{b_{1,2}}(\xi_2)$:
\begin{align*}
\langle \U_{b_{1,2}}(\xi_1)&\circ_1\U_{b_{1,2}}(\xi_2)|f\otimes
g\otimes h\rangle=
\langle \U_{b_{1,2}}(\xi_1)|\langle \U_{b_{1,2}}(\xi_2)|f\otimes
g\rangle\otimes h\rangle\\ &=\langle
\U_{b_{1,2}}(\xi_1)|2\,\xi_2^{i_3i_4}(\partial_{i_3}f)(\partial_{i_4}g)
\otimes h\rangle\\
&=4\,\xi_1^{i_1i_2}\partial_{i_1}(\xi_2^{i_3i_4}(\partial_{i_3}f)
(\partial_{i_4}g))(\partial_{i_2}h)
\\
&=4\,\xi_1^{i_1i_2}(\partial_{i_1}\xi_2^{i_3i_4})(\partial_{i_3}f)
(\partial_{i_4}g)(\partial_{i_2}h)
+4\,\xi_1^{i_1i_2}\xi_2^{i_3i_4}(\partial_{i_1}\partial_{i_3}f)
(\partial_{i_4}g)(\partial_{i_2}h)\\
&\phantom{mm}+4\,\xi_1^{i_1i_2}\xi_2^{i_3i_4}(\partial_{i_3}f)
(\partial_{i_1}\partial_{i_4}g)(\partial_{i_2}h)
\end{align*}

Now we turn to the computation of $\U_{b_{1,2}\circ_1 b_{1,2}}$. 
We have already computed that

\vbox{\hskip -1.5 cm
\includegraphics{fig-001.eps} 
}
 If we denote by $\Phi_1,\Phi_2$ and $\Phi_3$ the three graphs on the right 
hand side of the above equation, then
\begin{align*}
\langle\U_{\Phi_1}(\xi_1\otimes\xi_2)|f\otimes g\otimes h\rangle&=
\begin{xy}
,(0,0)*\cir(10,10){},
,(-7,-7)*{\bullet},(7,-7)*{\bullet}
,(1,-3);(-7,-7)**\crv{(0,-4)} ?>*\dir{>}
,(1,-3);(0,-10)*{\bullet}**\crv{(0,-4)} ?>*\dir{>}
,(3,3);(-7,-7)**\crv{(2,0)} ?>*\dir{>}
,(3,3);(7,-7)**\crv{(5,-5)} ?>*\dir{>}
,(4,5)*{\scriptstyle{\xi_1}}
,(2.5,-2)*{\scriptstyle{\xi_2}}
,(0,0.6)*{\scriptscriptstyle{i_1}}
,(5.2,1)*{\scriptscriptstyle{i_2}}
,(-3,-6.5)*{\scriptscriptstyle{i_3}}
,(1.5,-6)*{\scriptscriptstyle{i_4}}
,(-8,-9)*{\scriptstyle{f}}
,(0,-12)*{\scriptstyle{g}}
,(8,-9)*{\scriptstyle{h}}
\end{xy}\\
&=4\,
\xi_1^{i_1i_2}\xi_2^{i_3i_4}(\partial_{i_1}\partial_{i_3}f)
(\partial_{i_4}g)(\partial_{i_2}h)
\end{align*}

\begin{align*}
\langle\U_{\Phi_2}(\xi_1\otimes\xi_2)|f\otimes g\otimes h\rangle&=
\begin{xy}
,(0,0)*\cir(10,10){},
,(-7,-7)*{\bullet},(7,-7)*{\bullet}
,(1,-3);(-7,-7)**\crv{(0,-4)} ?>*\dir{>}
,(1,-3);(0,-10)*{\bullet}**\crv{(0,-4)} ?>*\dir{>}
,(3,3);(1,-3)**\crv{(3,0)} ?>*\dir{>}
,(3,3);(7,-7)**\crv{(5,-5)} ?>*\dir{>}
,(3.8,5)*{\scriptstyle{\xi_1}}
,(2.5,-3.8)*{\scriptstyle{\xi_2}}
,(1.8,1.4)*{\scriptscriptstyle{i_1}}
,(-3,-4)*{\scriptscriptstyle{i_3}}
,(5.9,-2.5)*{\scriptscriptstyle{i_2}}
,(-1,-6.5)*{\scriptscriptstyle{i_4}}
,(-8,-9)*{\scriptstyle{f}}
,(0,-12)*{\scriptstyle{g}}
,(8,-9)*{\scriptstyle{h}}
\end{xy}\\
&=4\,
\xi_1^{i_1i_2}(\partial_{i_1}\xi_2^{i_3i_4})(\partial_{i_3}f)
(\partial_{i_4}g)(\partial_{i_2}h)
\end{align*}
and
\begin{align*}
\langle\U_{\Phi_3}(\xi_1\otimes\xi_2)|f\otimes g\otimes h\rangle&=
\begin{xy}
,(0,0)*\cir(10,10){},
,(-7,-7)*{\bullet},(7,-7)*{\bullet}
,(-3,0);(-7,-7)**\crv{(-4,-4)} ?>*\dir{>}
,(-3,0);(0,-10)*{\bullet}**\crv{(-3,-4)} ?>*\dir{>}
,(3,3);(0,-10)**\crv{(3,0)} ?>*\dir{>}
,(3,3);(7,-7)**\crv{(5,-5)} ?>*\dir{>}
,(4.3,4.5)*{\scriptstyle{\xi_1}}
,(-3,1.2)*{\scriptstyle{\xi_2}}
,(1.1,-1.2)*{\scriptscriptstyle{i_1}}
,(-5.8,-3.8)*{\scriptscriptstyle{i_3}}
,(5.5,-1.3)*{\scriptscriptstyle{i_2}}
,(-1,-4.1)*{\scriptscriptstyle{i_4}}
,(-8,-9)*{\scriptstyle{f}}
,(0,-12)*{\scriptstyle{g}}
,(8,-9)*{\scriptstyle{h}}
\end{xy}\\
&=4\,
\xi_1^{i_1i_2}\xi_2^{i_3i_4}(\partial_{i_3}f)
(\partial_{i_4}\partial_{i_1}g)(\partial_{i_2}h)
\end{align*}
Therefore, $\langle \U_{b_{1,2}\circ_1 b_{1,2}}(\xi_1\otimes \xi_2)| f\otimes
g\otimes h\rangle=\langle \U_{b_{1,2}}\circ_1\U_{b_{1,2}}|f\otimes g\otimes
h\rangle$, as expected.
\end{example}

\subsection{$\U$ as a morphism of complexes}
\begin{lem}
$d_{1}(\U_\Gamma)=\U_{d_1\Gamma}$ and $d_{2}(\U_\Gamma)=\U_{d_2\Gamma}$. As a
consequence, $d_{CE}(\U_\Gamma)=\U_{d\Gamma}$
\end{lem}
\begin{proof}
The differential $d_1$ is defined as an adjoint both in the
Chevalley-Eilen\-berg complex and in the complex of graphs, and we have
already 
shown that $\U$ is a Lie algebra map ${\mathscr G}^{\,\bullet,\bullet}\to
CE^{\bullet,\bullet}(T^\bullet_{\rm poly},D^\bullet_{\rm poly})$. Then, the
compatibility of $\U$ with the differential $d_1$ immediately
follows from the fact that
$\U$ maps the graph $b_{0,2}$ (inducing $d_1$ on the complex of graphs) to the
multiplication $\mu$ on ${\cal C}^\infty({\setR}^d)$ (inducing
$d_1$ on the Chevalley-Eilenberg complex).
\par
The compatibility of $\U$ with the differential $d_2$ is an
immediate consequence of the finer
statement
$d_{{2,i}}(\U_\Gamma)=\U_{d_{2,i}\Gamma}$, where $d_{2,i}\colon
{\rm Hom}( (T^{\bullet}_{\rm
poly})^{\otimes q},D^\bullet_{\rm poly}[p])\to {\rm
Hom}((T^{\bullet}_{\rm
poly})^{\otimes (q+1)},D^\bullet_{\rm
poly}[p]) $ is defined as
\[
(d_{2,i}{\mathcal F})(\gamma_1\otimes\cdots
\otimes\gamma_i\otimes\gamma_{i+1}\otimes\cdots\otimes \gamma_{q+1})=
{\mathcal
F}(\gamma_1\otimes\cdots
\otimes(\gamma_i\bullet\gamma_{i+1})\otimes\cdots\otimes
\gamma_{q+1})
\]
and $d_{2,i}\colon{\mathscr G}^{\,p,q}\to{\mathscr G}^{\,p,q+1}$ is
the operator which
inserts an edge
 into the $i^{\rm th}$ vertex of a graph (see also \cite{Arnal,AM}).
Therefore, in order to prove the compatibility of $\U$ with the 
differential
$d_2$, we have to show that
 \begin{equation}\label{eq:bullet}
\U_{d_{2,i}\Gamma}
(\gamma_1\otimes\dots\otimes\gamma_{i}\otimes
\gamma_{i+1}\otimes\cdots\otimes\gamma_{q+1})=
\U_{\Gamma}(\gamma_1\otimes\dots\otimes(\gamma_{i}
\bullet\gamma_{i+1})\otimes
\cdots\otimes\gamma_{q+1}).
\end{equation}
This directly follows from the
definition of the Nijenhuis-Richardson pre-Lie operation, when
$\gamma_{i},\gamma_{i+1}$ are two homogeneous polyvector fields. 
Indeed, if 
$b_{1,m}$ denotes the $m$-corolla (i.e., the admissible graph with
one internal vertex, $m$ boundary vertices and $m$ edges), then the
definition of the Nijenhuis-Richardson
$\bullet$ operation is equivalent to the following identity:
\[
\U_{d_{2,1}b_{1,m}}(\gamma_1\otimes\gamma_2)=\U_{b_{1,m}}(\gamma_1\bullet\gamma_2).
\]
By linearity,
equation (\ref{eq:bullet}) holds for arbitrary (i.e., possibly
non-homogeneous) polyvector fields.
\end{proof}
\begin{rem}
Note that, if
$\deg(\gamma_{j})=k_j$,
then the evaluation of the terms in the left hand side of equation (\ref{eq:bullet})
corresponding to the different possible insertions of an
oriented edge in the vertex $i$ is non-trivial only
 when there are exactly $k^{}_i+1$ edges stemming from the vertex $i$ 
and $k^{}_{i+1}+1$ edges stemming from the vertex $i+1$.
\end{rem}

\begin{example} Let $b_{1,2}$ be the $2$-corolla.
The only nonzero homogeneous component of
$\U_{b_{1,2}}:T^\bullet_{\rm poly}\to D^\bullet_{\rm poly}$ is:
$$\U_{b_{1,2}}:T^{1}_{poly}\to D^1_{\rm poly},
$$ which is the natural
identification of 2-vector fields with 2-derivations.
The equations we want to verify are
\[
(d_{2,1}\U_{b_{1,2}})(\xi_1\otimes
\xi_2)=\U_{b_{1,2}}(\xi_1\bullet\xi_2)
=\U_{d_{2,1}b_{1,2}}(\xi_1\otimes\xi_2).
\]
These equations are trivial for homogeneous polyvector fields
unless 
$\deg(\xi_1\bullet\xi_2)=1$, i.e. unless $\deg\xi_1+\deg\xi_2=1$.
In these cases the first equation reduces to
\[
(d_{2,1}\U_{b_{1,2}})(\xi_1\otimes
\xi_2)=
\xi_1\bullet\xi_2
\]
Since the degree of a poly-vector field is at least $-1$, 
we have only two possibilities (up to a permutation),
namely $\deg\xi_1=1,\,\deg\xi_2=0$ or
$\deg\xi_1=2,\,\deg\xi_2=-1$.
Let us compute $\xi_1\bullet\xi_2$ in the first case. We find
\[
\xi_1\bullet\xi_2=
2\,\xi_1^{i_1i_2}(\partial_{i_1}\xi^{i_3}_2)
\partial_{i_2}\wedge\partial_{i_3}.
\]
Therefore,
\begin{align*}
\langle(d_{2,1}\U_{b_{1,2}})(\xi_1\otimes\xi_2)|f\otimes
g\rangle&=2\,\xi_1^{i_1i_2}(\partial_{i_1}\xi^{i_3}_2)
(\partial_{i_2}f)(\partial_{i_3}g)\\
&\phantom{mm}-2\,
\xi_1^{i_1i_2}(\partial_{i_1}\xi^{i_3}_2)
(\partial_{i_3}f)(\partial_{i_2}g)
\end{align*}
Now, we compute
$\langle\U_{d_{2,1}b_{1,2}}(\xi_1\otimes\xi_2)|f\otimes g\rangle$. We
have already computed

\vbox{\hskip -1.5 cm
\includegraphics{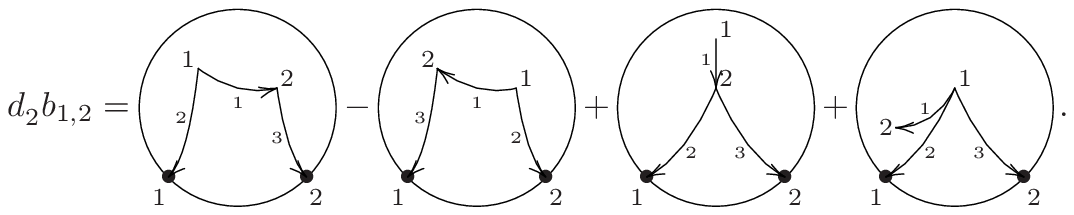}
}
If we denote the 4 graphs on the right hand side by
$\Phi_1,\dots,\Phi_4$, then the only non-zero pairings are
the following:

\begin{align*}
\langle\U_{\Phi_1}(\xi_1\otimes\xi_2)|f\otimes
g\rangle&=
\begin{xy}
,(0,0)*\cir(10,10){},
,(-7,-7)*{\bullet},(7,-7)*{\bullet}
,(-4,4);(-7,-7)**\crv{(-5,-5)} ?>*\dir{>}
,(-4,4);(4,2)**\crv{(0,1)} ?>*\dir{>}
,(4,2);(7,-7)**\crv{(5,-5)} ?>*\dir{>}
,(5,3.5)*{\scriptstyle{\xi_2}}
,(-5,5.3)*{\scriptstyle{\xi_1}}
,(0,0.6)*{\scriptscriptstyle{i_1}}
,(-6,-1)*{\scriptscriptstyle{i_2}}
,(4,-3)*{\scriptscriptstyle{i_3}}
,(-8,-9)*{\scriptstyle{f}}
,(8,-9)*{\scriptstyle{g}}
\end{xy}\\
&=2\,
\xi_1^{i_1i_2}(\partial_{i_1}\xi_2^{i_3})
(\partial_{i_2}f)(\partial_{i_3}g),
\end{align*}
and
\begin{align*}
\langle\U_{\Phi_2}(\xi_1\otimes\xi_2)|f\otimes
g\rangle&=
\begin{xy}
,(0,0)*\cir(10,10){},
,(-7,-7)*{\bullet},(7,-7)*{\bullet}
,(-4,4);(-7,-7)**\crv{(-5,-5)} ?>*\dir{>}
,(4,2);(-4,4)**\crv{(0,1)} ?>*\dir{>}
,(4,2);(7,-7)**\crv{(5,-5)} ?>*\dir{>}
,(5.4,3.7)*{\scriptstyle{\xi_1}}
,(-5,5.2)*{\scriptstyle{\xi_2}}
,(0,0.6)*{\scriptscriptstyle{i_1}}
,(-5.9,-1)*{\scriptscriptstyle{i_3}}
,(4,-3)*{\scriptscriptstyle{i_2}}
,(-8,-9)*{\scriptstyle{f}}
,(8,-9)*{\scriptstyle{g}}
\end{xy}\\
&=2\,
\xi_1^{i_1i_2}(\partial_{i_1}\xi_2^{i_3})
(\partial_{i_3}f)(\partial_{i_2}g)
\end{align*}

Therefore, if $\deg\xi_1=1$ and $\deg\xi_2=0$, we find
\begin{align*}
\langle{\U_{d_{2,1}b_{1,2}}}(\xi_1\otimes\xi_2)|f\otimes g\rangle&=
\langle{\U_{\Phi_1}}(\xi_1\otimes\xi_2)|f\otimes g\rangle
-\langle{\U_{\Phi_2}}(\xi_1\otimes\xi_2)|f\otimes g\rangle\\
&=2\,\xi_1^{i_1i_2}(\partial_{i_1}\xi_2^{i_3})
(\partial_{i_2}f)(\partial_{i_3}g)\\
&\phantom{mmm}-2\,\xi_1^{i_1i_2}(\partial_{i_1}\xi_2^{i_3})
(\partial_{i_3}f)(\partial_{i_2}g)\\
&=
\langle(d_{2,1}\U_{b_{1,2}})(\xi_1\otimes\xi_2)|f\otimes
g\rangle
\end{align*}
Now we consider the second possibility, namely $\deg\xi_1=2$
and $\deg\xi_2=-1$. In this case the Nijenhuis-Richardson bracket is
\begin{align*}
\xi_1\bullet\xi_2
&=3\,\xi_1^{i_1i_2i_3}(\partial_{i_1}\xi_2)
\partial_{i_2}\wedge\partial_{i_3}
\end{align*}
Therefore,
\begin{align*}
\langle(d_{2,1}\U_{b_{1,2}})(\xi_1\otimes\xi_2)|f\otimes
g\rangle=6\,
\xi_1^{i_1i_2i_3}(\partial_{i_1}\xi_2)
(\partial_{i_2}f)(\partial_{i_3}g)
\end{align*}
In computing
$\langle\U_{d_{2,1}b_{1,2}}(\xi_1\otimes\xi_2)|f\otimes g\rangle$
the only non trivial pairing is
\begin{align*}
{\U}_{\Phi_4}(\xi_1\otimes\xi_2)&=
\begin{xy}
,(0,0)*\cir(10,10){},
,(-7,-7)*{\bullet},(7,-7)*{\bullet}
,(0,2);(-7,-7)**\crv{(-3,-5)} ?>*\dir{>}
,(0,2);(-6,-2)**\crv{(-2,-2)} ?>*\dir{>}
,(0,2);(7,-7)**\crv{(3,-5)} ?>*\dir{>}
,(-7,-2)*{\scriptstyle{\xi_2}}
,(-3,0)*{\scriptscriptstyle{i_1}}
,(0,4)*{\scriptstyle{\xi_1}}
,(-2,-4.5)*{\scriptscriptstyle{i_2}}
,(2.5,-4.5)*{\scriptscriptstyle{i_3}}
,(-8,-9)*{\scriptstyle{f}}
,(8,-9)*{\scriptstyle{g}}
\end{xy}\\ \\
&=6\,
\xi_1^{i_1i_2i_3}(\partial_{i_1}\xi_2)(\partial_{i_2}f)
(\partial_{i_3}g)\\
&=d_{2,1}(\U_{b_{1,2}})(\xi_1\otimes\xi_2)|f\otimes
g\rangle
\end{align*}
Therefore, if $\deg\xi_1=2$ and $\deg\xi_2=-1$, we find
\begin{align*}
\langle{\U_{d_{2,1}b_{1,2}}}(\xi_1\otimes\xi_2)|f\otimes g\rangle=
\langle{\U_{\Phi_4}}(\xi_1\otimes\xi_2)|f\otimes g\rangle
=
\langle(d_{2,1}\U_{b_{1,2}})(\xi_1\otimes\xi_2)|f\otimes
g\rangle
\end{align*}
\end{example}

\section{The Maurer-Cartan equation on graphs and the
cocycle equation}\label{S:mccocycle}

We have seen in Section \ref{S:mclinf} that an $L_\infty$-morphism between
$T^\bullet_{\rm poly}$ and $D^\bullet_{\rm poly}$ can be viewed as a solution
of the vertical Maurer-Cartan equation 
on $CE^\bullet(T^\bullet_{\rm poly},D^\bullet_{\rm poly})$. Therefore,  
the differential graded Lie algebra morphism $\U$ reduces
the problem of writing an $L_\infty$-morphism between
$T^\bullet_{\rm poly}$ and $D^\bullet_{\rm poly}$ 
to
the problem of solving the vertical Maurer-Cartan equation on
the complex of graphs. This will in turn be equivalent to
solving the \emph{cocycle equation} which is the dual of the
vertical Maurer-Cartan equation.

\subsection{The dual dg-Lie-coalgebra structure}

We begin by passing from the graded Lie algebra of graphs
${\mathscr G}^{\,\bullet}$ to the dual Lie-coalgebra $({\mathscr
G}^{\,\bullet})^*$, where ${}^*$ means ``graded dual''. Since we
want to interpret also
elements in the Lie-coalgebra as graphs, we introduce a symmetric
non-degenerate pairing  
\[
\langle\,|\,\rangle\colon {\mathscr G}^{\,\bullet,\bullet}\otimes
{\mathscr G}^{\,\bullet,\bullet}\to {\mathbb C}
\]
 by declaring the set ${\mathcal G}^{\bullet,\bullet}$
of isomorphism classes of admissible graphs to be an orthonormal
basis of
${\mathscr G}^{\,\bullet,\bullet}$. This means that we are
identifying $({\mathscr G}^{p,q})^*$ with ${\mathscr
G}^{\,p,q}$. The pairing
$\langle\,|\,\rangle$ induces a differential graded
Lie coalgebra structure on
${\mathscr G}^{\,\bullet,\bullet}$ which we will now describe
explicitly. 
The \emph{codifferential} $d_2^*\colon{\mathscr G}^{\,p,q}\to
{\mathscr G}^{\,p,q-1}$ is defined by the equation
\[
\langle d_2^*\Phi|\Gamma\rangle=\langle
\Phi|d_2^{}\Gamma\rangle
\]
for any $\Phi,\Gamma$. Therefore
\begin{align*}
d_2^*\Phi&=\sum_\Gamma\langle d_2^*\Phi|\Gamma\rangle\Gamma
=\sum_\Gamma\langle \Phi|d_2^{}\Gamma\rangle\Gamma\\
&=\sum_\Gamma\sum_{\substack{e\hookrightarrow\Psi\twoheadrightarrow\Gamma\\
\text{internal}}}\pm\langle \Phi|\Psi\rangle\Gamma
=\sum_{\substack{ e\hookrightarrow\Phi\twoheadrightarrow\Gamma\\
\text{internal}}}\pm\Gamma\\
&=\sum_{\substack{e\in {\rm Edges}(\Phi)\\ \text{internal}}}\pm\Phi/e
\end{align*}
i.e., $d_2^*\Phi$ is the sum over all the internal edge contractions of $\Phi$.
Note that this is the Kontsevich graph homology differential
 \cite{Kfd}, p.109 (\cite{Kncg}, p.175). 
To determine an explicit expression
for the Lie co-bracket 
\[
\delta\colon {\mathscr
G}^{\,p,q}\to
\bigoplus_{{\substack{p_1+p_2=p\\q_1+q_2=q}}}{\mathscr
G}^{\,p_1,q_1}\otimes{\mathscr G}^{\,p_2,q_2}
\] 
we argue in a similar way:
\[
\bigl\langle \delta\Phi\bigr\rvert \Gamma\otimes\Psi\bigr\rangle=
\bigl\langle \Phi
\bigr\rvert [\Gamma,\Psi]\bigr\rangle
\]
so that
\begin{align*}
\delta\Phi&=\sum_{\Gamma,\Psi}\langle \delta\Phi|\Gamma\otimes
\Psi\rangle\Gamma\otimes \Psi
 =\sum_{\Gamma,\Psi}\bigl\langle \Phi
\bigr\rvert [\Gamma,\Psi]\bigr\rangle\Gamma\otimes \Psi\\
&=\sum_{\Gamma,\Psi}\sum_{\substack{\Gamma\hookrightarrow\Xi
\twoheadrightarrow\Psi\\
\text{boundary}}}\pm\langle
\Phi\vert\Xi\rangle
\bigl(\Gamma\otimes\Psi
\pm
\Psi\otimes\Gamma\bigr)\\
&=\sum_{\substack{\Gamma\normal \Phi\\
\text{boundary}}}\pm
\bigl(\Gamma\otimes\Phi/\Gamma
\pm
\Phi/\Gamma\otimes\Gamma\bigr).
\end{align*}
Note that the Lie cobracket $\delta$ coincides with the
anti-symmetrization of the
coproduct $\Delta$ introduced by Ionescu in \cite{pqft} (see also
\cite{cfg,linf}). 

\subsection{The Maurer-Cartan equation on graphs and the 
cocycle equation}\label{S:ce}

In this section we will start our study of the solutions of 
the vertical Maurer-Cartan equation
\[
d_2^{}\gamma+\frac{1}{2}[\gamma,\gamma]=0
\]
in the Hochschild-Kontsevich differential Lie algebra of graphs. 

Thanks to the non-degenerate pairing $\langle\,|\,\rangle$,
the above equation is equivalent to
\[
\bigl\langle
d_2^{}\gamma+\frac{1}{2}[\gamma,\gamma]\bigr\rvert
\Phi\bigr\rangle=0,\qquad\forall \Phi\in {\mathcal G}^2
\]
Now, we use the fact that ${\mathcal G}^1$ is a linear basis
of 
${\mathscr G}^{\,1}$ to write
\[
\gamma=\sum_{\G\in{\mathcal G}^1}W^{}_\Gamma\Gamma
\]
for some constants $W^{}_\Gamma$. The Maurer-Cartan equation then becomes
\[
\sum_\Gamma W^{}_\Gamma\bigl\langle
d_2^{}\Gamma\bigr\rvert
\Phi\bigr\rangle+\frac{1}{2}\sum_{\Gamma_1,\Gamma_2}
W^{}_{\Gamma_1}W^{}_{\Gamma_2}\bigl\langle
[\Gamma_1,\Gamma_2]\bigr\rvert
\Phi\bigr\rangle=0,\qquad\forall \Phi\in {\mathcal G}^2
\]
i.e.
\[
\sum_\Gamma W^{}_\Gamma\bigl\langle
\Gamma\bigr\rvert
d_2^*\Phi\bigr\rangle+\frac{1}{2}\sum_{\Gamma_1,\Gamma_2}
W^{}_{\Gamma_1}W^{}_{\Gamma_2}\bigl\langle
\Gamma_1\otimes\Gamma_2\bigr\rvert
\delta\Phi\bigr\rangle=0,\qquad\forall \Phi\in {\mathcal G}^2
\]
The map $W\colon\Gamma\mapsto W_\Gamma$ can be seen as a linear functional on
${\mathscr G}^{\,1}$, i.e., as an element of the space of
``weights''. Having this in mind, the above equation can
be rewritten as the
\emph{cocycle equation}
\begin{equation}\label{E:cceq}
W(d_2^*\Phi)+\frac{1}{2}W^{\otimes 2}(\delta\Phi)=0,\qquad\forall
\Phi\in{\mathcal G}^2
\end{equation}
which explicitly reads
\begin{equation}\label{E:cceq2}
\sum_{\substack{e\in E_\Phi\\
\text{internal}}}\pm W_{\Phi/e}+\sum_{\substack{\Gamma\normal \Phi\\
\text{boundary}}}\pm W_\Gamma W_{\Phi/\Gamma}
=0,\qquad\forall \Phi\in {\mathcal G}^2
\end{equation}
We will call any solution $W$
of the above equation a \emph{cocycle}.
The reason for this is that $W$ is actually a cocycle for the
cobar differential corresponding to the coalgebra structure
on 
${\mathscr G}^{\,\bullet,\bullet}$; see details in subsection
\ref{ss:cobar}. 
\begin{rem} In order to make direct contact with the original proof
by Kontsevich in \cite{K1}, we introduced the somehow artificial
self-duality of the space of graphs by imposing that graphs form an
orthonormal basis. Actually, the natural dual of the space of graphs
is the space of weights, and the cocycle equation should be
directly read as the dual Maurer-Cartan equation for the
dg-Lie-coalgebra structure on the space of weights.
\end{rem}
 Summing up what we have proved
so far in the paper, we obtain:
\begin{prop}\label{P:linf}
Let $W$ be any cocycle such that $W_{b_{0,2}}=1$. Then $\mathcal
F^\mu_{}=\break\sum_{\Gamma\in {\mathcal G}^1}W_\Gamma\U_\Gamma$ is
an
$L_\infty$-morphism between
$T^\bullet_{\rm poly}$ and $D^\bullet_{\rm poly}$. 
\end{prop}
\begin{rem} \label{rem:K-star-prod}
As seen in Remark \ref{R:defq}, 
we obtain a deformation quantization
of Poisson structures out of any ${\mathcal F}^\mu_{}$ such that
${\mathcal F}^\mu_{1}\colon T^m_{\rm poly}\to D^m_{\rm poly}$ is the
canonical identification of polyvector fields with polydifferential
operators. 
Since the restriction of $\sum_{\Gamma\in {\mathcal G}^1}
W_\Gamma\U_\Gamma$ to ${\rm Hom}(T^{m-1}_{\rm poly},D^{m-1}_{\rm
poly})$ is $W_{b_{1,m}}\U_{b_{1,m}}$ and
$\U_{b_{1,m}}(\xi)=m!\xi$, we get a star-product quantization 
of Poisson structures out of any cocycle with 
$W_{b_{1,m}}=\frac{1}{m!}$. The formula for the star-product is
\[
f\star_\alpha g=\sum_{\Gamma\in{\mathcal G}^1}
W_\Gamma\U_\Gamma(\exp\{\hbar\alpha\})(f,g)
=\sum_{q=0}^\infty\sum_{\Gamma\in{\mathcal G}^{1-q,q}}
\frac{\hbar^q}{q!}W_\Gamma\U_\Gamma(\alpha^{\wedge q}_{})(f,g).
\] 
Recall from Section \ref{S:gc} that 
${\mathcal G}^{1-q,q}$ is the
set of isomorphism classes of all admissible graphs $\Gamma$
with $q$ internal vertices and such that
$e(\Gamma)-v(\Gamma)+2=q$. Here
$v(\Gamma)$ and $e(\Gamma)$ are the number of vertices and
edges of $\Gamma$, respectively. So, if we denote by $n$
the number of internal vertices and by $m$ the number of
boundary vertices of $\Gamma$, a graph in ${\mathcal
G}^{1-q,q}$ satisfies $n=q$ and $e(\Gamma)=2n+m-2$.
Since $\alpha$ is a 2-vector field, the
admissible graphs actually appearing in the star-product formula must
have exactly two outgoing edges for any internal vertex, and
we get the identity 
\[
2n=e(\Gamma)=2n+m-2
\]
which implies $m=2$, as it should be since $\U_\Gamma$ must be an
element of $D^1_{\rm poly}$. 
Then, if we denote by ${\mathcal G}^1_{n,2}$
the set of isomorphism classes of degree one admissible graphs with $n$ 
internal and two boundary vertices, 
and for any graph $\Gamma$ in ${\mathcal G}_{n,2}^{1}$ we write
$B_{\Gamma,\alpha}$ for $\U_\Gamma(\alpha^{\wedge n})$,  the
formula for the star-product gets the familiar aspect (see
\cite{K1})
\[
f\star_\alpha g=\sum_{n=0}^\infty
\frac{\hbar^n}{n!}
\sum_{\Gamma\in{\mathcal G}_{n,2}^{1}}
W_\Gamma B_{\Gamma,\alpha}(f,g).
\]
\end{rem}

The problem of (local) deformation quantization is thus
reduced to solving the cocycle equation on graphs, i.e., to
determining a system of weights $W_\Gamma$ satisfying 
equation (\ref{E:cceq2}).
Such a system of weights has first been  found in
\cite{K1}, where the $W_\Gamma$ are defined as integrals on
suitable configuration spaces of points in a disk, naturally
associated to admissible graphs. For a  sigma model/quantum
field theory interpretation of these configuration space
integrals, we refer the reader to
\cite{CF,pqft}. An attempt towards a combinatorial solution
of the cocycle equation can be found in \cite{comb}.

\subsection{Connes-Kreimer convolution and Kontsevich
solution} The cocycle equation can be seen as a Connes-Kreimer
convolution with some vanishing coefficients. Namely,
equation (\ref{E:cceq2}) can be rewritten as
\begin{equation}\label{E:ckeq}
\begin{cases}
\displaystyle{\sum_{\substack{\Gamma\text{ proper}\\
\text{subgraph of }\Phi}}\!\!\!\! W_\Gamma
W_{\Phi/\Gamma} =0,\qquad\forall \Phi\in {\mathcal G}^2\,}
\\ \\
W_\Gamma W_{\Phi/\Gamma}=0, \qquad\text{unless
$\Gamma$ is an internal  edge or}\\
\phantom{W_\Gamma W_{\Phi/\Gamma}=0,}\qquad\text{ $\Gamma$
is a boundary normal subraph.}
\end{cases}
\end{equation}
The left-hand side of the first equation in (\ref{E:ckeq}) is 
Connes-Kreimer convolution $W\star W$, and the equation
$W\star W=0$ should be thought as an abstract version of the
Stokes formula on configuration spaces. Indeed, if 
${\rm Conf}$ denotes the configuration space functor mapping
an admissible graph
$\Gamma$ to the configuration space of the vertices of
$\Gamma$ into a manifold with boundary (in the interior or on
the boundary of the manifold depending on the type of the
vertices),
$\omega_\Gamma^{}$ denotes an angle form induced by an
embedding of $\Gamma$, and ${\setH}$ denotes the
Poincar\'e disk, then one can define the weights as $W=\langle
\overline{\rm Conf}(\setH)|\omega\rangle$, i.e.,
$
W_\Gamma=\int_{\overline{\rm Conf}_\Gamma(\setH)}
\omega_\Gamma^{}.
$
Then, 
\begin{multline*}
W\star W=\langle \overline{\rm Conf}(\setH)\otimes
\overline{\rm Conf}(\setH)|\omega\wedge
\omega\rangle\\=
\langle \partial \overline{\rm
Conf}(\setH)|\omega\rangle=\langle  \overline{\rm
Conf}(\setH)|d\omega\rangle=0,
\end{multline*}
where the last equality follows from the fact that angle
forms are closed: $d\omega=0$. This shows that the first
equation in (\ref{E:ckeq}) is  automatically satisfied as soon as
one defines the weights as configuration space integrals. The
proof of the cocycle equation is then reduced to the checking
of the vanishing of the configuration space integrals
corresponding to the terms in the second line of
(\ref{E:ckeq}), as in \cite{K1}. See the Appendix for
details on configuration space integrals in Kontsevich's proof.

\subsection{Configuration spaces and the cobar
construction}\label{ss:cobar} 
We now explain in which sense a solution $W$ to the cocycle
equation (\ref{E:cceq2}) is a cocycle. Recall \cite{Quillen}
that a DGLA structure on a vector space ${\mathfrak g}$ can
be seen as a codifferential $Q$ which is a coderivation on
$C({\mathfrak g}[1])=\bigoplus_{n\geq 1}S^n({\mathfrak
g}[1])$, the cofree coassociative cocommutative coalgebra
without counit cogenerated by ${\mathfrak g}[1]$. In our case,
${\mathfrak g}={\mathscr G}^{\,\bullet}$ and the
codifferential
$Q$ is
$d_2^*+\delta$. We have remarked that the cobraket $\delta$
is actually the antisymmetrization of a coproduct $\Delta$,
which is dual to the pre-Lie operation
$\circ$. This implies \cite{pqft} that the differential
$d_2^*+\delta$ on $C({\mathscr G}^{\,\bullet}[1])$ is
induced by the differential $D=d_2^*+\Delta$ on $T_+({\mathscr
G}^{\,\bullet}[1])=\bigoplus_{n\geq 1}T^n({\mathscr
G}^{\,\bullet}[1])$, the cofree coassociative
coalgebra without counit cogenerated by
${\mathscr G}^{\,\bullet}[1]$. Note that
$T_+({\mathscr G}^{\,\bullet}[1])$ is nothing else
than the graded vector space underlying the cobar complex
$F({\mathscr G}^{\,\bullet})$ of
${\mathscr G}^{\,\bullet}$ (see
\cite{J}, p.111). The differential $D$ can be seen as a
perturbation of the differential $d_2^*$ of the cobar complex.
Note that $D=d_2^*+\Delta$ is the differential 
yielding the cohomology of Feynman-Kontsevich graphs 
\cite{pqft,cfg,linf}.
The system of weights $W$ is as linear functional on
$T_+({\mathscr G}^{\,\bullet}[1])$; from this point of
view, the cocycle equation is nothing but the equation
${}^tDW=0$, where
${}^tD$ is the dual differential on ${\rm Hom}(T_+({\mathscr
G}^{\,\bullet}[1]);{\setC})$. In other words, $W$ is
a ${}^tD$-cocycle.
\par

The operator $D$ provides a nice interpretation of the fact
that Kontsevich's weights do indeed provide a cocycle. Recall
from the previous subsection 
that
$W_\Gamma$ is defined using a rule
$\omega$
which maps
an admissible graph $\Gamma$  
to a differential form $\omega_\Gamma$ on the configuration
space of the vertices of $\Gamma$ in the Poincar\'e disk:
\[
\omega\colon {\mathscr
G}^{\,\bullet}\to\Omega^{\bullet}\bigl(\overline{\rm
Conf}(\setH);\setR\bigr)
\]
Such a form is the angle form; 
for details, in addition to \cite{K1}, 
see \cite{Po}. 
The map $\omega$ can be extended to a map between the cobar
constructions on graphs and on the space of differential
forms on configuration spaces. Using the familiar fact that
differential forms on a product are the tensor product of
differential forms on the factors, we get a map
\[
\omega\colon \left\{T_+({\mathscr
G}^{\,\bullet}[1]),D\right\}\to
\left\{
\Omega^\bullet\left(T_+\bigl(\overline{\rm
Conf}(\setH)\bigr);\setR\right),d\right\},
\]
which is a map of complexes.\footnote{This is, in a slightly
different formulation from the original paper, Theorem 3.2
of \cite{pqft}.}  The weights $W$ on the cobar complex of
graphs are then defined as $W=\langle T_+\bigl(\overline{\rm
Conf}(\setH)\bigr)|\omega\rangle$, where the angle form
corresponding to a tensor product of admissible graphs is
integrated over the corresponding component in
$T_+\bigl(\overline{\rm Conf}(\setH)\bigr) $.
The proof 
that Kontsevich coefficients yield a cocycle then reads as
follows: for any admissible graph $\Gamma$,
\[
\langle {}^tD W|\Gamma\rangle=\langle W|D\Gamma\rangle=
\langle T_+\bigl(\overline{\rm
Conf}(\setH)\bigr)|\omega_{D\Gamma}^{} \rangle=
\langle T_+\bigl(\overline{\rm
Conf}(\setH)\bigr) |d\omega_{\Gamma}^{} \rangle=0,
\]
since angle forms are closed.

\section{The tree-level approximation}
In the main body of the paper we have reviewed how an $L_\infty$-morphism
between $T^\bullet_{\rm poly}$ and $D^\bullet_{\rm poly}$
can be written as a sum over admissible graphs,
leading to a deformation quantization of Poisson structures. We
now show how the DGLA point of view we adopted provides a nice
algebraic setting to discuss the tree-level approximation of
Kontsevich's star-product. Namely, we discuss the possibility of
writing a star-product formula as a sum over the family of Kontsevich's
forests, i.e. admissible graphs whose underlying non-oriented graphs
become acyclic\footnote{A non-oriented graph is called acyclic when all of
its connected components are simply connected.} when one
removes the boundary vertices. There are indeed several hints
that such a formula may exist. For instance, constant Poisson
structures admit the Moyal's star-product quantization which can be seen
as a sum over graphs with no edge ending on an internal
vertex, and Michael Polyak has shown in
\cite{Po} that the Kontsevich's star-product for linear Poisson
structures is gauge-equivalent to the Gutt's star-product,
which can be seen as a sum over the family of Kontsevich's
forests.
\footnote{The reader interested in Kontsevich's star-product 
for linear Poisson structures is also addressed to \cite{ABM} and \cite{D}.}
It should be remarked that Polyak's proof relies on the
linearity of the coefficients of the Poisson structures and so
does not provide an argument to decide whether the family
of forests could be universal, in the sense that a formula
written as a sum over forests could work for any Poisson
structure. Additional evidence towards the existence of such a
universal formula is provided in
\cite{comb}, where the cocycle condition is investigated from a
combinatorial point of view (see also \cite{IS}).
\par
In this final section we sketch how an analysis of the
tree-level approximation can be obtained by looking
at the cohomology of the complex of admissible graphs.
The first thing we note is that Kontsevich's forests span a
subcomplex $({\mathscr F}^{\,\bullet},d_2^{})$ of $({\mathscr
G}^{\,\bullet},d_2^{})$, since inserting an internal
edge does not change the homotopy type of a graph.
On the other hand, ${\mathscr F}^{\,\bullet}$
is not a DGLA,
since the Lie bracket of two forest graphs in not necessarly
a forest, as the following example shows
\[
\left[
\begin{xy}
,(0,0)*\cir(10,10){},
,(-7,-7)*{\bullet},(7,-7)*{\bullet}
,(0,2);(-7,-7)**\crv{(-3,-5)} ?>*\dir{>}
,(2,-2);(-7,-7)**\crv{(-3,-6)} ?>*\dir{>}
,(0,2);(2,-2)**\crv{(.5,0)} ?>*\dir{>}
,(2,-2);(7,-7)**\crv{(3,-4)} ?>*\dir{>}
,(0,3.3)*{\scriptstyle{1}}
,(-4,-2.5)*{\scriptscriptstyle{1}}
,(1.5,0.7)*{\scriptscriptstyle{2}}
,(-1.7,-6)*{\scriptscriptstyle{3}}
,(3.2,-5)*{\scriptscriptstyle{4}}
,(3.1,-1.5)*{\scriptstyle{2}}
,(-8,-9)*{\scriptstyle{1}}
,(8,-9)*{\scriptstyle{2}}
\end{xy}\,\,,\,\,\begin{xy}
,(0,0)*\cir(10,10){},
,(-7,-7)*{\bullet},(7,-7)*{\bullet}
,(0,2);(-7,-7)**\crv{(-3,-5)} ?>*\dir{>}
,(0,2);(7,-7)**\crv{(3,-5)} ?>*\dir{>}
,(1,3)*{\scriptstyle{1}}
,(-4,-2.5)*{\scriptscriptstyle{1}}
,(4,-2.5)*{\scriptscriptstyle{2}}
,(-8,-9)*{\scriptstyle{1}}
,(8,-9)*{\scriptstyle{2}}
\end{xy}\right]\,\,
=\,\,
\begin{xy}
,(0,0)*\cir(10,10){},
,(-7,-7)*{\bullet},(7,-7)*{\bullet}
,(0,-4);(-7,-7)**\crv{(-1,-5)} ?>*\dir{>}
,(0,-4);(0,-10)*{\bullet}**\crv{(-1,-5)} ?>*\dir{>}
,(3,3);(0,-4)**\crv{(3,0)} ?>*\dir{>}
,(4.5,-2.5);(7,-7)**\crv{(5.3,-5)} ?>*\dir{>}
,(3,3);(4.5,-2.5)**\crv{(3.5,0)} ?>*\dir{>}
,(4.5,-2.5);(0,-4)**\crv{(3,-2)} ?>*\dir{>}
,(4,4)*{\scriptstyle{1}}
,(-1,-3)*{\scriptstyle{3}}
,(1.8,0.6)*{\scriptscriptstyle{1}}
,(4.6,0.6)*{\scriptscriptstyle{2}}
,(3.2,-3.4)*{\scriptscriptstyle{3}}
,(4.7,-5.1)*{\scriptscriptstyle{4}}
,(-3.8,-4.7)*{\scriptscriptstyle{5}}
,(5.6,-2.2)*{\scriptstyle{2}}
,(.5,-6.4)*{\scriptscriptstyle{6}}
,(-8,-9)*{\scriptstyle{1}}
,(0,-12)*{\scriptstyle{2}}
,(8,-9)*{\scriptstyle{3}}
\end{xy}
\,\,
+\cdots
\]
Therefore we cannot immediately reduce the tree-level
approximation to a cocycle equation for forests. But one can
show (see for instance
\cite{Kom}) that the inclusion
$\iota\colon{\mathscr F}^{\,\bullet}
\hookrightarrow {\mathscr G}^{\,\bullet}$ is a
quasi-isomorphism of complexes. Then, by homotopical transfer
of structure (see
\cite{KS,Markl})
${\mathscr F}^{\,\bullet}$ inherits from the
DGLA structure of admissible graphs an $L_\infty$-structure
making ${\mathscr F}^{\,\bullet}$ and
${\mathscr G}^{\,\bullet}$ homotopy equivalent  
$L_\infty$-algebras, via an $L_\infty$-morphism
$\iota_\infty\colon {\mathscr F}^{\,\bullet}\to{\mathscr
G}^{\,\bullet}$. The homotopy equivalence $\iota_\infty$
is a non-linear map,\footnote{Note that, since ${\mathscr
F}^{\,\bullet}$ is not closed under the Lie bracket, the
quadratic component $\iota_{\infty,2}$ is non-trivial.} whose linear term
coincides with
the inclusion $\iota$.  
By Kontsevich-Schles\-singer deformation
theory \cite{K1} the solutions of the Maurer-Cartan equation  
on the DGLA ${\mathscr
G}^{\,\bullet}$ bijectively correspond (up to gauge
equivalence) to solutions of the Maurer-Cartan equation 
$
\sum_{n=1}^\infty[\gamma^{\wedge n}]^{}_n/n!=0
$
on the
$L_\infty$-algebra ${\mathscr F}^{\,\bullet}$.
In particular, Kontsevich's solution on ${\mathscr
G}^{\,\bullet}$ will induce a solution
$
\gamma_\infty^{}=\sum_{T\in {\mathcal F}^1_{}}W_{\infty,T}T,
$
which is a sum over forest graphs. Note that, due to the
non-linearity of the homotopical isomophism $\iota_\infty^{}$,
the weights $W_\infty^{}$ are not simply the restriction of the
weights $W$ to forest graphs. To relate the sum over forests
$\gamma_\infty^{}$ to the tree-level approximation of
Kontsevich's star-product, we compose the homotopical
isomorphism $\iota_\infty^{}$ with the DGLA map $\U$, to get 
an $L_\infty$-morphism $\U_\infty\colon{\mathscr
F}^{\,\bullet}\to CE^\bullet(T^\bullet_{\rm
poly},D^\bullet_{\rm poly})$, whose linear part
$\U_{\infty,1}$ is just the restriction of $\U$ to forest
graphs. By construction, evaluating
$\U_\infty^{}$ on $\gamma_\infty^{}$ is the same thing as
evaluating $\U$ on the Kontsevich's solution, so we have found
an expression for the formality $L_\infty$-morphism between
$T^\bullet_{\rm poly}$ and $D^\bullet_{\rm poly}$ which is a sum
over forest graphs. Explicitly,
\begin{align*}
\U_\infty^{}(\gamma_\infty^{})&=\sum_{n=1}^\infty\frac{1}{n!}
\U_{\infty,n}(\gamma_\infty^{\wedge
n})\\&=\sum_{T\in {\mathcal
F}^1_{}}W_{\infty,T}\U_T+\frac{1}{2}\sum_{T_1,T_2\in
{\mathcal
F}^1_{}}W_{\infty,T_1}W_{\infty,T_2}\U_{\infty,2,T_1,T_2}
+\cdots
\end{align*}
Evaluating this expression on $\exp\{\hbar\alpha\}$, for a
given Poisson bivector field
$\alpha$, one finds
\begin{thm} The Kontsevich's star-product formula quantizing
a Poisson structure $\alpha$ on $\setR^d$ can be rewritten as
\begin{align*}
f\star_\alpha^{}g&=\underbrace{\sum_{n=0}^\infty
\frac{\hbar^n}{n!}
\sum_{T\in{\mathcal F}_{n,2}^1}
W_{\infty,T}B_{T,\alpha}(f,g)}_{\text{tree level
approximation}}+\\
&\phantom{mmm}+\underbrace{\sum_{k=2}^\infty\frac{1}{k!}\sum_{n=0}^\infty
\hbar^n\sum_{T_i\in{\mathcal F}^1_{m_i,n_i}}
\left(\prod_{i=1}^k\frac{
W_{\infty,T_i}}{n_i!}\right)B_{\infty,T_1,\dots,T_k,\alpha}
(f,g)}_{\text{correction terms coming from loops}},
\end{align*}
where the sum
ranges over the $n_i$ and $m_j$ such that $n_1+\cdots+
n_k=n$ and
$m_1+\cdots+m_k=2+k-1$.
\end{thm}
In more colloquial terms, the theorem above states that it is
possible to write a star-product formula as a
sum over forest graphs; this formula is given
by a Kontsevich-type formula restricted to forest graphs (but
with different weights), the so-called tree-level or
semi-classical approximation, plus corrections
corresponding to contributions coming from non
simply-connected graphs in the original Kontsevich's formula.
In the sum-over-forests reformulation, these contributions
are written in terms of the multilinear operations
$T_1\wedge\cdots \wedge
T_k\mapsto B_{\infty,T_1,T_2,\ldots,T_k,\alpha}$ on forest
graphs.
\par
 It should be
remarked that, since the tree-level approximation depends on
the choice of the homotopy equivalence $\iota_\infty^{}$,
it is far from being unique. Thus we have several different
tree-level approximations of Kontsevich's formula, which
correspond to different systems of weights $W_\infty$; clearly
the various tree-level approximations are gauge-equivalent to
each other. Recasted into these terms, Polyak's result
mentioned at the beginning of this section states that for a
linear Poisson structure $\alpha$ there exist a choice of
$\iota_\infty^{}$ (i.e., a choice of the weights $W_\infty$)
such that the tree-level approximation is exact. The problem
of deciding whether for an arbitrary Poisson structure
$\alpha$ there exist an exact tree-level approximation
remains open. Yet, it is not unlikely that a more refined
study of the $L_\infty$-algebra structure on ${\mathscr
F}^{\,\bullet}$, of the multilinear operations
$\iota_{\infty,k}\colon
\wedge^k{\mathscr F}^{\,\bullet}\to {\mathscr
G}^{\,\bullet}[1-k]$ and, consequently, of the multilinear
operations $\U_{\infty,k}$ and $B_{\infty,k,\alpha}$ on
forest graphs may turn to be useful to get a better
understanding of this problem. Also it would be interesting
to have an interpretation of the weights $W_\infty^{}$ in
terms of the geometry of suitable configuration spaces.

\section*{Appendix: Configuration spaces and Kontsevich's
cocycle}\label{S:config}
\setcounter{section}{1}
\renewcommand{\thesection}{\Alph{section}}
\addcontentsline{toc}{section}{\thesection\hspace{1em}
Appendix: Configuration spaces and Kontsevich's
cocycle}

 In this section we review
Kontsevich's solution of the cocycle equation, yielding the
weights $W_\Gamma^{}$ in the graph expansion of the formality $L_\infty$-morphism.
We will not be interested into the physical interpretation behind this
construction
(see \cite{CF,pqft} for details); rather we will stress the r\^ole
of the dg-Lie-coalgebra structure and its relations with Stokes'
theorem for configuration spaces of points in a disk, as a crucial
ingredient of Kontsevich's proof.

\subsection{Configuration spaces and differential forms associated with
graphs}\label{ss:conf-sp} Let
${\mathcal G}^{}_{n,m}$ be the set
of (isomorphism classes of) admissible graphs with $n$
internal  and
$m$ boundary vertices. With any
$\G\in{\mathcal G}^{}_{n,m}$, Kontsevich associates a manifold with
corners
$\overline{C}_{n,m}$ and a closed differential form $\omega^{}_\Gamma$ on
$\overline{C}_{n,m}$. As the notation suggests, the manifold
$\overline{C}_{n,m}$ depends only on the number and type of vertices of $\Gamma$.
Indeed, the manifold $\overline{C}_{n,m}$ is the Fulton-MacPherson
compactification of the configuration space of $n+m$ points in the disk
\[
\Delta=\{z\in\setC\text{ such that } ||z||\leq1\},
\] 
with $n$ points in the interior of the disk and $m$ points on the
boundary. Under the usual
identification of the boundary of the disk with
$\setR\cup\infty$, one can assume that the boundary points are represented by
$0=t_1<t_2<\cdots<t_m=1$. This means that we are actually fixing an
additional point ``$\infty$'' on the boundary of the disk. Since the group
of complex automorphisms of $\Delta$ fixing $\infty$ has real dimension
$2$, the real dimension of $\OC_{n,m}$ is $2n+m-2$.
The differential form
$\omega^{}_\Gamma$ is defined by the following construction: a point of
$\overline{C}_{n,m}$ determines a configuration of points in $\Delta$. Use
this configuration to draw a copy of
$\Gamma$ in the disk, drawing each edge of $\Gamma$ as a geodesic segment (in
the usual hyperbolic metric of $\Delta$). Then each edge $e$ of $\Gamma$
determines an angle $\phi_e$, which is the angle between the edge $e$ and
the geodesic joining the starting point of $e$ with the point $\infty$ on
the
boundary of the disk. As usual, $\phi_e$ is locally well defined up to an
integral multiple of $2\pi$ so that the 1-form $d\phi_e$ is well defined. The
closed form $\omega^{}_\Gamma$ is then defined as the normalized wedge
product of all the 1-forms $d\phi_e$ with $e$ ranging in the set $E_\Gamma$
of edges of $\Gamma$:
\[
\omega^{}_\Gamma=\frac{1}{(2\pi)^{|E_\G|}}
\bigwedge_{e\in E_\G}d\phi_e.
\]
\par
If $\Gamma$ has bidegree $(p,q)$ 
then $|E_\Gamma|=2n+m-(p+q)-1=\dim\OC_{n,m}+1-(p+q)$
and the differential form $\omega^{}_\Gamma$ is a top-dimensional differential
form on $\OC_{n,m}$ precisely when $\Gamma$ has total degree
1 (see Remark \ref{rem:K-star-prod}).
 Having a manifold and a top-dimensional differential form
on it, we can get a real number out of them: define the
 weight $W_\Gamma$, for an admissible graph of total degree
1, as
\[
W_\Gamma=\int_{\OC_{n,m}}\omega^{}_\Gamma.
\]
On the other hand, if $\Gamma$ has total degree 2, then the
number of edges of $\Gamma$ is $2n+m-3$, so that we get a
top-form on $\OC_{n,m}$ by taking the differential
$d\omega^{}_\Gamma$, in this case. Stokes' theorem holds
 for manifold with corners so we get, for any
admissible graphs with $n$ internal and $m$ boundary 
vertices having total degree 2, the identity
\[
     \int\limits_{\p \OC_{n,m}}\omega^{}_\Gamma
     =\int\limits_{\OC_{n,m}}d\omega^{}_\Gamma
     =0\,\,\,.
\]
In the next section we will relate the boundary integral on
the left to the weights $W_\G$ of degree 1 admissible
graphs. 
To do this we need a description of the codimension one boundary
components of $\OC_{n,m}$.

\subsection{Normal subgraphs and the boundary behavior of the
$\omega^{}_\Gamma$'s}\label{ss:boundary}

Boundary components of
$\OC_{n,m}$ describe the various ways in which a subset of
points in a given configuration can  collapse to a single
point. Since we are dealing with configuration
spaces of a manifold with boundary, there are two kinds of
degenerations (and so two kinds of boundary components).
Degenerations of the first kind are those in which a set $S$
consisting of
$n_1\geq 2$ internal points collapse to an internal point.
In this case, the collapsing makes up a
``bubble'' and in the limit we are left with a copy of
${\mathbb P}^1(\setC)$ with $n_1$ points on it, 
joined to a disk with $n-n_1$ internal points left on it. 
The joining point of the sphere and the disk is an
additional marked point for both spaces. Moreover, the direction pointing
towards $\infty$ is a distinguished direction in the disk, and leads to a
distinguished tangent vector at $\infty$ on ${\mathbb
P}^1(\setC)$. 
So we get the boundary component
$$\p_S\OC_{n,m}\simeq \OC_{n_1}\times \OC_{n-n_1+1,m}\ ,$$
where $\OC_{n_1}$ denotes the configuration space of $n_1$ points
distinct from $\infty$ on ${\mathbb P}^1(\setC)$, 
up to projective transformations fixing a tangent direction at
$\infty$. 
Note that, since the group of these projective transformations 
has real dimension $3$, the manifold $\OC_{n_1}$ has real dimension $2n_1-3$,
and 
\[
\dim(\p_S\OC_{n,m})=\dim(\OC_{n,m})-1.
\]
 Degenerations of the second kind are those in which a set
$S$ consisting of $n_1$ internal points 
and a non-empty set $S'$ consisting of
$m_1$ boundary points, with $2n_1+m_1\geq 2$ collapse to a boundary point. In
this case, the collapse produces a new disk joined to the original one by a
point in the boundary, which is the $\infty$ point in the new disk. This way
we get  the boundary component
$\p_{S,S'}\OC_{n,m}\simeq \OC_{n_1,m_1}\times \OC_{n-n_1,m-m_1+1}$. 
Note that also in this case we find
\[
\dim(\p_{S,S'}\OC_{n,m})=\dim(\OC_{n,m})-1.
\]
\begin{rem} Since it may be a source of confusion, we
explicitly note that $\OC_{n,0}$ and $\OC_n$ are different
spaces (even their dimensions are different). 
\end{rem}
Note that the above description of the boundary of $\OC_{n,m}$ produced an
asymmetry: there are factors in the boundary components which are associated
to admissible graphs (namely, the factors of the form $\OC_{n',m'}$) and
factors which are not associated to graphs (namely, the factors of the form
$\OC_{n'}$). 
To remedy this asymmetry, 
we add to admissible graphs a copy of
${\mathcal G}_{n,0}$, which we will denote by ${\mathcal G}_{n}$. 
The difference between graphs  in ${\mathcal G}_{n}$ and 
graphs in ${\mathcal G}_{n,0}$ is that the former have to be thought 
as embedded into ${\mathbb P}^1({\mathbb C})$, which has no boundary, 
whereas the latter have to be thought as embedded into $\Delta$, 
which has a boundary, but the graphs have no boundary vertices. 
To a graph $\Gamma$ in ${\mathcal G}_{n}$ 
we will associate the configuration space
$\OC_{n}$ and a differential form $\omega^{}_\Gamma$
on $\OC_{n}$ by the same rule as before.
Note that, in contrast with what happened for
$\OC_{n,m}$, the differential form
$\omega_\Gamma^{}$ will be top-dimensional on $\OC_{n}$
precisely when $\Gamma$ has total degree 2.

Each subgraph $\Phi$ of $\Gamma$ determines a subset of points, the
vertices of $\Phi$, denoted by $v(\Phi)$. 
Conversely, each subset $S$
(or $S,S'$) of the vertices of $\Gamma$
determines a full subgraph of $\Gamma$ by the rule ``edges of $\Phi$ are
those edges of $\Gamma$ which join the vertices in $S$ (or in
$S\cup S'$)''. The
subgraphs corresponding to the sets $S$ and $S,S'$ will be denoted $\Gamma_S$
and $\Gamma_{S,S'}$ respectively.
\par
The boundary behavior of the differential forms $\omega^{}_\Gamma$ can be
described in a clear way by means of the notion of graph extension
introduced in Section \ref{S:graphs}. 
Consider a graph $\Gamma$ in ${\mathcal G}_{n,m}$ with total degree $2$, 
and let $S$ be a subset of its internal vertices,  with $|S|\geq 2$.
Corresponding to $S$ we have a boundary component
\[
\iota_S\colon \OC_{|S|}\times \OC_{n-|S|+1,m}\hookrightarrow \OC_{n,m}.
\]
In addition, we have the subgraph $\Gamma_S$ of $\Gamma$ and the quotient
graph $\Gamma/\Gamma_S$. 
Assume that $\Gamma_S$ is a normal subgraph of $\Gamma$. 
Then $\Gamma_S\in{\mathcal G}_{|S|}$ and
$\Gamma/\Gamma_{S}\in {\mathcal G}_{n-|S|+1,m}$. 
Therefore we have a differential form
$\omega^{}_{\Gamma_S}\wedge \omega^{}_{\Gamma/\Gamma_S}$ on $\OC_{|S|}\times
\OC_{n-|S|+1,m}$. Since $\Gamma$
is obtained by internal composition of $\Gamma_S$ and
$\Gamma/\Gamma_{S}$, and the internal composition $\bullet$
has total degree $-1$, then
$\Gamma/\Gamma_S$ has degree 1 precisely when $\Gamma_S$ has
degree
$2$, and $\omega^{}_{\Gamma_S}\wedge 
\omega^{}_{\Gamma/\Gamma_S}$ is a top-dimensional
differential form on
$\OC_{|S|}\times
\OC_{n-|S|+1,m}$ in this case.
\par
Similarly, if $(S,S')$ are subsets of the internal
and boundary vertices of $\Gamma$ respectively, with $2|S|+|S'|\geq 2$, then
we have a boundary component
\[
\iota_{S,S'}\colon \OC_{|S|,|S'|}\times \OC_{n-|S|,m-|S'|+1}\hookrightarrow
\OC_{n,m}
\]
and, in case the subgraph $\Gamma_{S,S'}$ is normal, we also have the
differential form
$\omega^{}_{\Gamma_{S,S'}}\wedge\omega^{}_{\Gamma/\Gamma_{S,S'}}$ on
$\OC_{|S|,|S'|}\times \OC_{n-|S|,m-|S'|+1}$. Since the
boundary composition of graphs has total degree zero, then
$\Gamma_{S,S'}$ has degree 1 precisely when also
$\Gamma/\Gamma_{S,S'}$ has degree one. In this case the
differential form $\omega^{}_{\Gamma_{S,S'}}\wedge\omega^{}_{\Gamma/\Gamma_{S,S'}}$ on
$\OC_{|S|,|S'|}\times \OC_{n-|S|,m-|S'|+1}$ is
top-dimensional. 
One of the essential ingredients of
Kontsevich's proof are the following
\emph{compatibility equations}: let $\Gamma$ be a degree 2
admissible graph, an let $S$, $S'$ be such that
$\deg\Gamma/\Gamma_S=\deg\Gamma/\Gamma_{S,S'}=1$. 
Then 
\begin{equation*}
\begin{matrix}
\iota^*_S\omega^{}_\Gamma=\omega^{}_{\Gamma_S}\wedge\omega^{}_{\Gamma/\Gamma_S}\\
\\
\iota^*_{S,S'}\omega^{}_\Gamma=\omega^{}_{\Gamma_{S,S'}}\wedge
\omega^{}_{\Gamma/\Gamma_{S,S'}}
\end{matrix}
\end{equation*}
The above compatibility equations tell us that
$\omega$ can be seen as an Euler-Poincar\'e map with respect to short exact
sequences of admissible graphs \cite{pqft}.
\par
 It remains to be seen what happens when one considers a full subgraph of
$\Gamma$ which is not normal. Here is where the second main ingredient,
the \emph{Kontsevich vanishing theorem}, comes in: if $S$ or $(S,S')$ correspond to non-normal subgraphs, then
\begin{equation*}
\begin{matrix}
{\displaystyle{\int_{\OC_{|S|}\times
\OC_{n-|S|+1,m}}\iota^*_S\omega^{}_\Gamma=0}}
\\
{\displaystyle{\int_{\OC_{|S|,|S'|}\times
\OC_{n-|S|,m-|S'|+1}}\iota^*_{S,S'}\omega^{}_\Gamma=0}}
\end{matrix}
\end{equation*}
Indeed, non-normal subgraphs correspond to graphs with ``bad edges'' in
Kontsevich's original notations, see \cite{K1}, p.27 . 
There is another vanishing theorem we will need in the next section: 
if $S$ is a subset of
internal vertices of $\Gamma$ with $|S|\geq 3$ then the integral of
$\omega^{}_{\Gamma_S}$ on $\OC_{|S|}$ vanishes. As a consequence,
\begin{equation*}
\begin{matrix}
\\
{\displaystyle{\int_{\OC_{|S|}\times
\OC_{n-|S|+1,m}}\iota^*_S\omega^{}_\Gamma=0\qquad\text{ if } |S|\geq 3}}
\\
\end{matrix}
\end{equation*}
This means that, in order to have a non-vanishing integral
we have to take $|S|=2$. But then, the requirement
$\deg\Gamma_S^{}=2$ forces $\Gamma_S$ to be an edge with two
internal vertices; the configuration space $\OC_2$ is
diffeomorphic to $S^1$ and
\[
\int_{\OC_2}\omega_{\Gamma_S}
=\frac{1}{2\pi}\int_{S^1}d\phi_e=1
\]
\par
Summing up, from the compatibility and the
vanishing equations above we have the following behavior
for the boundary integrals of the
$\omega_\Gamma^{}$'s: let
$\Gamma$ be an admissible graph of total degree 2, with $n$
internal and
$m$ boundary vertices. Then
\begin{equation*}
\begin{matrix}
{\displaystyle{\int_{\OC_{|S|}\times
\OC_{n-|S|+1,m}}\iota^*_S\omega^{}_\Gamma=\begin{cases}
W_{\Gamma/\Gamma_S^{}}^{}\text{ if 
$\Gamma_S$ is an internal edge of $\Gamma$}\\
\phantom{mmmi}\text{ such
that
$\Gamma/\Gamma_S^{}$ is admissible}\\
\\0\phantom{m}\text{ elsewhere}\end{cases}}}
\end{matrix}
\end{equation*}
and
\begin{equation*}
\begin{matrix}
{\displaystyle{\int_{\OC_{|S|,|S'|}\times
\OC_{n-|S|,m-|S'|+1}}\iota^*_S\omega^{}_\Gamma=
\begin{cases}
W_{\Gamma_{S,S'}}^{}W_{\Gamma/\Gamma_{S,S'}}^{}
\text{ if $\Gamma_{S,S'}$ is a normal}\\
\phantom{mmmmmmiiii}\text{ subgraph of $\Gamma$ of total}\\
\phantom{mmmmmmiiii}\text{ degree 1}\\
\\0\phantom{m}\text{ elsewhere}\end{cases}}}
\end{matrix}
\end{equation*}

\subsection{Stokes' theorem and the cocycle equation}
We now turn back to Stokes' theorem:
\[
\int\limits_{\p \OC_{n,m}}\omega^{}_\Gamma=0,
\]
for any admissible graph $\Gamma$ of total degree 2.
According to the results of the previous subsection, this
equation can be rewritten as 
\begin{align*}
0&=\sum_{S}  \int\limits_{\p_S \OC_{n,m}}\omega^{}_\Gamma
      + \sum_{S, S'}  \int\limits_{\p_{S,S'}
       \OC_{n,m}}\omega^{}_\Gamma \\
&=\sum_{e\in
E_\Gamma}\pm W_{\Gamma/e}+
\sum_{\substack{S, S'\\
\Gamma_{S,S'}\normal\Gamma}}\pm W_{\Gamma_{S,S'}}W_{\Gamma/\Gamma_{S,S'}}
\end{align*}
where the $\pm$ signs depend on the orientation of the boundary components of
$\OC_{n,m}$. Due
to the bijective correspondence between subsets of the vertices of 
$\Gamma$ and full
subgraphs of $\Gamma$, we have
\[
\sum_{\substack{S, S'\\
\Gamma_{S,S'}\normal\Gamma}}\pm W_{\Gamma_{S,S'}}W_{\Gamma/\Gamma_{S,S'}}=
\sum_{\substack{\Phi\text{ full, boundary}\\
\Phi\normal\Gamma}}\pm W_{\Phi}W_{\Gamma/\Phi}
\]
Moreover, we have seen is Section \ref{S:graphs} that a normal
 subgraph is necessarily full; so, summing up, we
have found the equation
\begin{equation*}
\sum_{e\in
E_\Gamma}\pm W_{\Gamma/e}+\sum_{\substack{\Phi\text{boundary}\\
\Phi\normal\Gamma}}\pm W_{\Phi}W_{\Gamma/\Phi}=0,
\end{equation*}
that is,
\[
W(d_2^*\Gamma)+\frac{1}{2}W^{\otimes
2}(\delta_{}\Gamma)=0,\qquad\forall
\Gamma\in {\mathcal G}^2,
\] 
i.e. the Kontsevich integral $W$ is a cocycle. 
\begin{rem}
One explicitly computes $W_{b_{0,2}}=1$ and
$W_{b_{1,m}}=1/m!$. Thus the element
$\sum_{\Gamma\in{\mathcal G}^1}W_\Gamma\,\Gamma$ is a
solution of the Maurer-Cartan equation on graphs fulfilling
the conditions from Proposition \ref{P:linf} and Remark
\ref{rem:K-star-prod}.  In particular Kontsevich's weights define
a star-product deformation quantization of Poisson structures.
\end{rem}



\def\cprime{$'$}

\end{document}